\documentclass[10pt,legno]{article}
\usepackage{amsmath, amssymb}
\usepackage{latexsym}
\usepackage{amssymb,epsfig,color}
\usepackage{graphicx,epsfig}
\usepackage[toc,page]{appendix}
\begin{document}
\newcommand{\qed}{\hfill \ensuremath{\square}}
\newtheorem{thm}{Theorem}[section]
\newtheorem{cor}[thm]{Corollary}
\newtheorem{lem}[thm]{Lemma}
\newtheorem{prop}[thm]{Proposition}
\newtheorem{defn}[thm]{Definition}
\newcommand{\proof}{\vspace{1ex}\noindent{\em Proof}. \ }
\newtheorem{pro}[thm]{proof}
\newtheorem{ide}[thm]{Idee}
\newtheorem{rem}[thm]{Remark}
\newtheorem{ex}[thm]{Example}
\bibliographystyle{plain}
\numberwithin{equation}{section}
\numberwithin{equation}{section}
\newcounter{saveeqn}
\newcommand{\subeqn}{\setcounter{saveeqn}{\value{equation}}%
 \stepcounter{saveeqn}\setcounter{equation}{0}
\renewcommand{\theequation}{\mbox{\arabic{section}.\arabic{saveeqn}\alph{=
equation}}}} 

\newcommand{\reseteqn}{\setcounter{equation}{\value{saveeqn}}%
\renewcommand{\theequation}{\arabic{section}.\arabic{equation}}}

\def\nm{\noalign{\medskip}}
\def\u{{\mathbf u}}
\def\e{{\mathbf e}}
\def\T{{\mathbf T}}
\def\aa{{\mathbf a}}
\def\A{{\mathbf A}}
\def\B{{\mathbf B}}
\def\v{{\mathbf v}}
\def\g{{\mathbf g}}
\def\I{{\mathbf I}}
\def\f{{\mathbf f}}
\def\n{{\mathbf n}}
\def\T{{\mathbf T}}
\def\N{{\mathbf N}}
\def\M{{\mathbf M}}
\def\w{{\mathbf w}}
\def\f{{\mathbf f}}
\def\x{{\mathbf x}}
\def\F{{\mathbf F}}
\def\GG{{\mathbf G}}
\def\H{{\mathbf H}}
\def\U{{\mathbf U}}
\def\y{{\mathbf y}}
\def\t{{\mathbf t}}
\def\T{{\mathbf T}}
\def\bvarphi{\boldsymbol{\varphi}}
\def\bpsi{\boldsymbol{\psi}}
\def\bphi{\boldsymbol{\phi}}
\def\ta{\boldsymbol{\tau}}
\def\G{{\mathbf \Gamma}}
\def\etaa{{\boldsymbol \eta}}
\newcommand{\Om}{\Omega}
\newcommand{\om}{\omega}
\newcommand{\Real}{\mathbb{R}}
\newcommand{\nuu}{\tilde{{\nu}}}
\newcommand{\xe}{\tilde{x}}
\newcommand{\ye}{\tilde{{y}}}
\newcommand{\bohm}{{\partial}{\ohm}}
\newcommand{\la}{\langle}
\newcommand{\ra}{\rangle}
\newcommand{\ms}{\mathcal{S}_\ohm}
\newcommand{\mk}{\mathcal{K}_\ohm}
\newcommand{\mks}{\mathcal{K}_\ohm ^{\ast}}
\newcommand{\grad}{\bigtriangledown}
\newcommand{\ds}{\displaystyle}
\newcommand{\pf}{\medskip \noindent {\sl Proof}. ~ }
\newcommand{\p}{\partial}
\renewcommand{\a}{\alpha}
\newcommand{\z}{\zeta}
\newcommand\q{\quad}
\newcommand{\pd}[2]{\frac {\p #1}{\p #2}}
\newcommand{\pdl}[2]{\frac {\p^2 #1}{\p #2}}
\newcommand{\dbar}{\overline \p}
\newcommand{\eqnref}[1]{(\ref {#1})}
\newcommand{\na}{\nabla}
\newcommand{\ep}{\epsilon}
\newcommand{\vp}{\varphi}
\newcommand{\fo}{\forall}
\newcommand{\Scal}{\mathcal{S}}
\newcommand{\BScal}{\boldsymbol{\mathcal{S}}}
\newcommand{\BDcal}{\boldsymbol{\mathcal{D}}}
\newcommand{\Dcal}{\mathcal{D}}
\newcommand{\Kcal}{\mathcal{K}}
\newcommand{\BKcal}{\boldsymbol{\mathcal{K}}}
\newcommand{\K}{\boldsymbol{\mathbb{K}}}
\newcommand{\Ecal}{\mathcal{E}}
\newcommand{\Ncal}{\mathcal{N}}
\newcommand{\BNcal}{\boldsymbol{\mathcal{N}}}
\newcommand{\Abar}{\overline A}
\newcommand{\Rcal}{\mathcal{R}}
\newcommand{\Lcal}{\mathcal {L}}
\newcommand{\Tcal}{\mathcal {T}}
\newcommand{\Gcal}{\mathcal {G}
}
\newcommand{\Cbar}{\overline C}
\newcommand{\Ebar}{\overline E}
\newcommand{\RR}{\mathbb{R}}
\newcommand{\CC}{\mathbb{C}}
\newcommand{\NN}{\mathbb{N}}
\newcommand{\Z}{\mathbb{Z}}

\title{An Asymptotic Expansion for  Perturbations in the Displacement Field   due to the Presence of  Thin Interfaces}

\date{}

\author{ Jihene Lagha \thanks{  Department of Mathematics, Faculty of Sciences, University of Tunis El Manar, Tunis 2092, Tunisia (lagha.jihene@yahoo.fr) } \and Habib Zribi \thanks{  Department of Mathematics, College of Sciences, University of Hafr Al Batin,
P.o. 1803, Hafr Al Batin 31991, Saudi Arabia (zribi.habib@yahoo.fr)}}

\maketitle

\begin{abstract} We derive an asymptotic expansion for two-dimensional displacement field associated to thin elastic inhomogeneities having no uniform thickness. Our derivation is rigorous and based on layer potential techniques. We extend these techniques to determine
a relationship between traction-displacement measurements and the  shape of the object.
\end{abstract}

\noindent {\footnotesize {\bf Mathematics subject classification
(MSC2000):} 35B30, 35C20, 31B10}


\section{Introduction and statement of main results}
Let $D$ be a bounded $\mathcal{C}^{2,\eta}$ domain in $\RR^2$ for some $\eta>0$. For a given $\ep > 0$,
let $D_\ep$ be an $\ep-$perturbation of $D, i.e.$,  there is $h\in \mathcal{C}^{1}(\p D)$ such that $\p D_\ep$
is given by
\begin{equation*}
\ds \p D_{\ep}:=\Big\{\tilde {x}:\tilde{x}=x+\ep h(x)\n(x), x\in \p D \Big\},
\end{equation*}
where $\n(x)$ is the outward normal to the domain $D$. We assume that $h(x)\geq C> 0$ for all $x \in \p D$.

Consider a homogeneous isotropic  elastic object occupying $\RR^2$. Suppose  that the thin layer $D_\ep \backslash \overline{D}$ lies inside $\RR^2$.
 Let the constants $(\lambda_0,\mu_0),(\lambda_1,\mu_1)$, and $(\lambda_2,\mu_2)$ denote the Lam\'e coefficients of $\RR^2\backslash \overline{D}_\ep$, $D$, and $D_\ep \backslash \overline{D}$, respectively. It is always assumed  that
\begin{align*}
\ds \mu_j >0,\q\q\lambda_j +\mu_j>0\q \mbox{ for } j=0,1,2.
\end{align*}
We also assume that
\begin{align*}
(\lambda_0-\lambda_j)(\mu_0-\mu_j)\geq 0,\q\q \Big((\lambda_0-\lambda_j)^2+(\mu_0-\mu_j)^2\neq 0\Big)\q \mbox{ for }j=1,2.
\end{align*}

Let $\mathbb{C}_0$,  $\mathbb{C}_1$, and  $\mathbb{C}_2$ be  the elasticity tensors for $\RR^2\backslash \overline {D}_\ep$, $D$, and $D_\ep \backslash \overline{D}$, respectively,  which are  given by
\begin{align}\label{tensor-C}
\ds \mathbb{C}_m=\lambda_m \I\otimes\I+2 \mu_m \mathbb{I}, \q\q m=0,1,2,
\end{align}
where $\mathbb{I}$ is the identity  $4$-tensor and $\I$ is the identity $2$-tensor (the $2 \times 2$ identity matrix).

We define
\begin{align*}
\ds \mathbb{C}_\ep:=\mathbb{C}_0\chi_{\RR^2\backslash \overline{D}_\ep}+\mathbb{C}_2\chi_{D_\ep \backslash \overline{D}}+\mathbb{C}_1\chi_{D},\q\q
\mathbb{C}:=\mathbb{C}_0\chi_{\RR^2\backslash \overline{D}}+\mathbb{C}_1\chi_{D},
\end{align*}
where $\chi_{D}$ is the indicator function of $D$.

Let $\u_\ep$ denote the displacement field in the presence of the thin $D_\ep \backslash \overline{D}$, that is,  the solution to
\begin{equation}\label{Main-Pb-ep}
\left\{
  \begin{array}{lll}
   \ds \nabla \cdot \big(\mathbb{C}_\ep \widehat{\nabla}\u_\ep \big)=0 & \mbox{ in } \RR^2,\\
    \nm \ds \u_\ep(x)-\H(x)=O(|x|^{-1}) & \mbox{ as }|x|\rightarrow \infty,
  \end{array}
\right.
\end{equation}
where  $\widehat{\nabla}\u_\ep=\frac{1}{2}\big(\nabla \u_\ep+\nabla \u_\ep^{T}\big)$ is the strain tensor and $\H$
is a vector-valued function satisfying $\nabla \cdot \big(\mathbb{C}_0 \widehat{\nabla}\H\big)=0$ in $\RR^2$. Here and throughout the paper $\M^T$
denotes the transpose of the matrix $\M$.

The elastostatic system corresponding to the Lam\'e constants
$(\lambda_j, \mu_j)$ is defined by
\begin{align}\label{elastostatic-system}
\ds \mathcal{L}_{\lambda_j,\mu_j}\w:= \mu_j\Delta \w+(\lambda_j+\mu_j) \nabla \nabla \cdot \w.
\end{align}
The corresponding conormal derivative $ \displaystyle{\pd{\w}{\nu_j}}$ associated with $(\lambda_j,\mu_j)$ is defined by
\begin{align}\label{conormal-derivative}
\ds \pd{\w}{\nu_j}:= \lambda_j(\nabla \cdot \w)\n +\mu_j \big(\nabla \w+\nabla \w ^{T}\big)\n=\lambda_j(\nabla \cdot \w)\n +2\mu_j \widehat{\nabla}\w \,\n.
\end{align}
The problem \eqref{Main-Pb-ep} is equivalent to the
following problem
\begin{equation}\label{equation-u}
\left\{
  \begin{array}{ll}
   \ds  \mathcal{L}_{\lambda_0,\mu_0}\u_\ep=0 \q  \mbox{  in }\RR^2\backslash \overline{D}_\ep,  \\
    \nm \ds \mathcal{L}_{\lambda_2,\mu_2}\u_\ep=0  \q   \mbox{  in }{D}_\ep\backslash \overline{D},  \\
     \nm \ds \mathcal{L}_{\lambda_1,\mu_1}\u_\ep=0  \q  \mbox{  in }{D},  \\
   \nm \ds  \u_\ep\big|_{-}=\u_\ep\big|_{+} \q \q \mbox{ on } \p D, \q\q\u_\ep\big|_{-}=\u_\ep\big|_{+}  \q \q\mbox{ on  } \p D_\ep,\\
   \nm\ds \pd{\u_\ep}{\nu_1}\Big|_{-}=\pd{\u_\ep}{{\nu_2}}\Big|_{+} \q \mbox{ on } \p D, \q\mbox{ } \pd{\u_\ep}{\nu_2}\Big|_{-}=\pd{\u_\ep}{{\nu_0}}\Big|_{+} \q \mbox{ on } \p D_\ep,\\
    \nm \ds \u_\ep(x)-\H(x)=O(|x|^{-1}) \q\q \mbox{ as }|x|\rightarrow \infty.
  \end{array}
\right.
\end{equation}
The notation  $ \u_\ep|_{\pm}$ on $\p D$ denote the limits from outside and inside of $D$, respectively. We will sometimes use $\u_\ep^e$ for $\u_\ep|_{+}$ and $\u_\ep^i$ for $\u_\ep|_{-}$.

The first achievement of this work,  a rigorous derivation of  the asymptotic expansion of $\u_{\ep}|_{\Om}$ as $\ep\rightarrow 0$ where $\Om$ is  a bounded region away from $\p D$, based on   layer
potential techniques.
\begin{thm} \label{Main-theorem} Let $\u_\ep$ be the solution to \eqref{equation-u}. Let $\Om$ be a bounded region away from $\p D$. For $x \in \Om$,  the following pointwise   asymptotic expansion  holds:
\begin{align}\label{Main-asymptotic}
\ds  \u_\ep(x)  =\u(x)+\ep \u_1(x)+o(\ep),
\end{align}
where the remainder $o(\ep)$ depends only on $(\lambda_j, \mu_j)$ for j=0,1,2,
the $\mathcal{C}^2$-norm of $X$, the $\mathcal{C}^1$-norm of $h$, and $dist (\Om, \p D)$,
 $\u$ is the unique solution to
\begin{equation}\label{Main-Pb}
\left\{
  \begin{array}{lll}
   \ds \nabla \cdot \big(\mathbb{C} \widehat{\nabla}\u\big)=0 & \mbox{ in } \RR^2,\\
    \nm \ds \u(x)-\H(x)=O(|x|^{-1}) & \mbox{ as }|x|\rightarrow \infty,
  \end{array}
\right.
\end{equation}
and  $\u_1$ is the unique solution of the following transmission problem:
\begin{equation}\label{equation-u-1}
\left\{
  \begin{array}{ll}
   \ds  \mathcal{L}_{\lambda_0,\mu_0}\u_{1}=0 & \mbox{  in } \RR^2\backslash \overline{D},  \\
    \nm \ds \mathcal{L}_{\lambda_1,\mu_1}\u_{1}=0& \mbox{  in }{D},  \\
   \nm \ds  \u_1\big|_{-}-\u_1\big|_{+}= h\Big((\mathbb{K}_{0,1}-\mathbb{K}_{2,1})\widehat{\nabla}\u^{i}\Big)\n& \mbox{ on  } \p D,\\
   \nm\ds \pd{\u_1}{\nu_1}\Big |_{-}-\pd{\u_1}{\nu_0}\Big|_{+}=
   \frac{\p }{\p \ta}\Big(h\big[(\mathbb{M}_{2,1}-\mathbb{M}_{0,1})\widehat{\nabla}\u^{i}\big]\ta\Big) & \mbox{ on  } \p D,  \\
   \nm \ds \u_1(x)=O(|x|^{-1}) & \mbox{ as }|x|\rightarrow \infty,
  \end{array}
\right.
\end{equation}
with  $\ta$ is the tangential vector to $\p D$,
\begin{align*}
\ds \mathbb{M}_{l,k}&:=\frac{\lambda_l (\lambda_k+2\mu_k)}{\lambda_l+2\mu_l} \I \otimes \I+2 \mu_k \mathbb{I}+ \frac{4(\mu_l-\mu_k)(\lambda_l+\mu_l)}{\lambda_l+2\mu_l} \I\otimes (\ta \otimes \ta),\\
\nm \ds \mathbb{K}_{l,k}&:=\frac{\mu_l(\lambda_k-\lambda_l)+2(\mu_l-\mu_k) (\lambda_l+\mu_l)}{\mu_l(\lambda_l+2\mu_l)} \I \otimes \I+2\big(\frac{\mu_k}{\mu_l}-1\big) \mathbb{I}\\
\nm\ds &\q+ \frac{2(\mu_k-\mu_l)(\lambda_l+\mu_l)}{\mu_l(\lambda_l+2\mu_l)} \I\otimes (\ta \otimes \ta).
\end{align*}
\end{thm}

By taking $\lambda_2=\lambda_1$ and $\mu_2=\mu_1$, we reduce our problem to that proposed in \cite{JFH}. So it is obvious
to obtain the asymptotic expansion of the displacement field  resulting from small perturbations of
the shape of an elastic inclusion
already derived in \cite[Theorem 1.1]{JFH}. Our asymptotic expansion is valid in
the cases of thin interfaces with  high contrast parameters and can be derived similarly to \cite{KZ2},
for more details on behaviors of the  leading and  first  order terms $\u$ and $\u_1$, we refer the reader to \cite[Chapter 2]{ABGKLW}.
In \cite{ BBFM,BF}
 the authors derive asymptotic expansions
for the  boundary displacement field in the case $\lambda_0=\lambda_1$ and $\mu_0=\mu_1$ in both cases
 of isotropic and anisotropic thin elastic inclusions,
 the approach they use, based on energy estimates, variational  approach,  and
 fine regularity estimates for solutions of elliptic systems obtained by Li and Nirenberg \cite{LN}.
 Unfortunately,  this method does not seem to work in our case.

Let $\v$ be the solution of the following problem:
\begin{equation}\label{v}
\left\{
  \begin{array}{ll}
   \ds  \mathcal{L}_{\lambda_0,\mu_0}\v=0 & \mbox{  in }\RR^2\backslash \overline{D},  \\
    \nm \ds \mathcal{L}_{\lambda_1,\mu_1}\v=0  &  \mbox{  in }{D},  \\
   \nm \ds  \v|_{-}=\v|_{+}   &\mbox{ on  } \p D, \\
   \nm\ds \pd{\v}{{\nu}_1}\Big|_{-}=\pd{\v}{{\nu_0}}\Big|_{+} &\mbox{ on  } \p D,  \\
    \nm \ds   \v (x)-\F(x)=O(|x|^{-1})&\mbox{ as } |x|\rightarrow \infty.
  \end{array}
\right.
\end{equation}
As a consequence of the  theorem \ref{Main-theorem}, we obtain the following relationship between
traction-displacement measurements and the deformation $h$.
\begin{thm} \label{second-theorem} Let  $\u_\ep$, $\u$, and $\v$ be the solutions
 to \eqref{equation-u}, \eqref{Main-Pb}, and  \eqref{v}, respectively. Let $S$
  be a Lipschitz closed curve enclosing $D$  away from $\p D$.
The following asymptotic expansion holds:
\begin{align}\label{asymptotic-traction-displacement}
\ds &\int_{S}\big(\u_\ep-\u\big)\cdot \pd{ \F}{\nu_0}d\sigma-\int_{S}\big(\pd{\u_\ep}{\nu_0}-\pd{ \u}{\nu_0}\big)\cdot \F d\sigma\nonumber\\
 \nm \ds & \q =\ep\int_{\p D}h\bigg(\big([\mathbb{M}_{0,1}-\mathbb{M}_{2,1}]
 \widehat{\nabla}\u^i\big)\ta\cdot \widehat{\nabla} \v^i \ta +\big([\mathbb{K}_{2,1}-\mathbb{K}_{0,1}]
  \widehat{\nabla}\u^i\big)\n\cdot (\mathbb{C}_1 \widehat{\nabla} \v^i )\n \bigg)d\sigma+o(\ep),
\end{align}
where the remainder $o(\ep)$ depends only   on $(\lambda_j, \mu_j)$ for j=0,1,2,
the $\mathcal{C}^2$-norm of $X$, the $\mathcal{C}^1$-norm of $h$, and $dist (S, \p D)$. The dot denotes the scalar product in $\RR^2.$
\end{thm}

The asymptotic expansion in \eqref{asymptotic-traction-displacement} can be
used to design algorithms to identify certain properties of its thin elastic like location and thickness
 based on traction-displacement  measurements. To
do this, we refer to asymptotic formulae related to  measurements
in the same spirit, far-field data, currents,  generalized polarization
tensors, elastic moment tensors,  multistatic
response  at single or multiple frequencies,  and modal measurements that have been obtained in a series
of recent papers \cite{AEEKL, AGKLS, AKLZ1, AKLZ2, KKL, LLZ, LY, Zribi1}.

Our techniques in this paper seem simpler to implement and
have the big advantage to derive high order terms in
 the asymptotic formulae and allow a
generalization to $3$-dimensional interface problems by using  \cite{KZ1, KZ2}.

This paper is organized as follows. In section 2, we review  some
 basic facts on the layer potentials of the Lam\'e system and  derive
a representation formula for the solution of the  problem \eqref{equation-u}. In section 3,
we derive asymptotic expansions of layer potentials. In section 3,
based on layer potentials techniques   we rigorously derive the  asymptotic expansion for  perturbations in the displacement field
 and find the  relationship between traction-displacement measurements
 and the deformation $h$ (Theorem \ref{Main-theorem} \& Theorem \ref{second-theorem}).
\section{Representation formula}
Let us  review some well-known properties of the layer potentials for the
elastostatics. The theory of layer potentials has
been developed in relation to boundary value problems in a Lipschitz domain.

Let
\begin{align*}
\ds \Psi:=\Big \{\bpsi: \p_i \bpsi_j+\p_j \bpsi_i=0,\q 1\leq i,j\leq 2\Big \}.
\end{align*}
or equivalently,
\begin{align*}
\ds \Psi=\mbox{span}\Bigg\{\theta_1(x):={1 \brack 0}, \theta_2(x):={0 \brack 1},\theta_3(x):={x_2 \brack -x_1}\Bigg\}.
\end{align*}
Introduce the space
\begin{align*}
\ds \displaystyle L^2_{\Psi}(\p D):=\Big \{\f \in L^2(\p D): \int_{\p D} \f\cdot \bpsi ~d\sigma=0 \mbox{ for all } \bpsi \in \Psi \Big\}.
\end{align*}
In particular, since $\Psi$ contains constant functions, we get
\begin{align*}
 \ds \int_{\p D} \f d\sigma=0
\end{align*}
for any $\f \in L^2_{\Psi}(\p D)$. The following fact is useful later.
\begin{align}\label{L-Psi}
\ds \mbox{If} \q \w \in W^{1,\frac{3}{2}}(D)  \q\mbox{satisfies}\q \Lcal_{\lambda_0,\mu_0} \w=0 \q \mbox{in } D,\q\mbox{then}\q \pd{\w}{\nu_0}\Big|_{\p D}\in L^2_{\Psi}(\p D).
\end{align}

The  Kelvin matrix of fundamental  solution $\G_j$ for the  Lam\'e  system $\Lcal_{\lambda_j,\mu_j}$ in $\RR^2$ is known to be
\begin{align}\label{Kelvin}
\ds \G_j(x)=\frac{A_j}{2\pi}\log |x|\I-\frac{B_j}{2\pi} \frac{x\otimes x }{|x|^2},\q x \neq 0,
\end{align}
where
\begin{align*}
\ds A_j=\frac{1}{2}\Big(\frac{1}{\mu_j}+\frac{1}{2\mu_j+\lambda_j}\Big) \q \mbox{and}\q B_j=\frac{1}{2}\Big(\frac{1}{\mu_j}-\frac{1}{2\mu_j+\lambda_j}\Big).
\end{align*}
The single and double layer potentials of the density function  $\bphi$ on $L^2(\p D)$
associated with the Lam\'e parameters $(\lambda_j, \mu_j)$  are defined by
 \begin{align}
\ds \BScal_{j,D} [\bphi](x)
 =& \int_{\partial D} \G_j(x-y) \bphi(y)d\sigma(y), \quad x \in
 \RR^2,\label{single-layer}\\
\nm \ds \BDcal_{j,D} [\bphi] (x)
=&\int_{\p D}\Bigg(\lambda_j \nabla_y \cdot \G_j(x-y)\otimes \n(y)\nonumber\\
\nm\ds &\q\q \q+ \mu_j \Big(\big[\nabla_y \G_j(x-y)\n(y)\big]^{T}+\nabla_y \G_j^{T}(x-y)\n(y)\Big)\Bigg)\bphi(y)d\sigma(y)\nonumber\\
\nm\ds :=&\int_{\p D}\mathbb{K}_j(x-y)\bphi(y)d\sigma(y), \quad x \in \RR^2 \setminus \p
D. \label{double-layer}
\end{align}
The followings are the well-known properties of  the single and double layer potentials due to Dahlberg, Keing, and Verchota \cite{DKV}. Let $D$ be a Lipschitz  bounded domain in $\RR^2$. Then we have
\begin{align}
\ds \pd{\BScal_{j,D} [\bphi] }{\nu_j}\Big |_{\pm}(x) & = \Big (\pm \frac{1}{2} \I
+ \BKcal_{j,D}^* \Big ) [\bphi] (x) \quad
\mbox{a.e. } x \in \p D \label{nuS}, \\
 \nm \ds \BDcal_{j,D} [\bphi] \big|_{\pm} (x) & = \Big(\mp \frac{1}{2}
\I + \BKcal_{j,D} \Big) [\bphi](x) \quad \mbox{a.e. } x \in \p D,
\label{doublejump-h}
\end{align}
where $\BKcal_{j,D}$ is defined  by
\begin{align*}
\ds \BKcal_{j,D} [\bphi](x)={p.v.}\int_{ \p D}\mathbb{K}_j(x-y) \bphi(y) d\sigma(y) \quad \mbox{a.e. } x \in \p D,
\end{align*}
 and  $\BKcal_{j,D}^*$ is the adjoint operator of $\BKcal_{j,D} $, that is,
\begin{align*}
\ds\BKcal_{j,D}^* [\bphi](x)={p.v.}\int_{ \p D}\mathbb{K}_j^{T}(x-y)\bphi(y) d\sigma(y) \quad \mbox{a.e. } x \in \p D.
\end{align*}
Here ${p.v.}$ denotes the Cauchy principal value.

Let $\BDcal_{j,D}^{\sharp}$  be the standard double layer potential which is  defined for any $\bphi \in L^2(\p D)$ by
\begin{align}\label{dcal-sharp}
\ds \BDcal_{j, D}^{\sharp}[\bphi](x)=\int_{\p D} \pd{\G_j(x-y)}{\n(y)}\bphi(y)d\sigma(y), \q x\in \RR^2\backslash {\p D}.
\end{align}
The following lemma holds, for more details see \cite{JFH}.

\begin{lem}\label{D-sharp} Let $D$ be a bounded Lipschitz domain in $\RR^2$. For $\bphi \in L^2(\p D)$
\begin{align}
 \ds \BDcal^{\sharp}_{j,D} [\bphi] \big|_{\pm} (x) & = \Big(\mp \frac{1}{2\mu_j}\I\pm  B_j \n \otimes \n
 + \BKcal^{\sharp}_{j,D} \Big) \bphi(x) \quad \mbox{a.e. } x \in \p D,\label{D-Sharp-n}\\
\nm \ds \pd{\BScal_{j,D} [\bphi]}{\n} \Big |_{\pm}(x) & = \Big (\pm \frac{1}{2\mu_j}\I \mp  B_j \n \otimes \n
+ \big(\BKcal_{j,D}^\sharp\big)^* \Big ) \bphi (x) \quad
\mbox{a.e. } x \in \p D,\label{S-Sharp-n}
\end{align}
where $\BKcal^\sharp_{j,D}$ is defined  by
\begin{align*}
\ds \BKcal^\sharp_{j,D} [\bphi](x)=\mbox{p.v.}\int_{ \p D}\pd{}{\n(y)} \G_j(x-y) \bphi(y) d\sigma(y) \quad \mbox{a.e. } x \in \p D,
\end{align*}
and  $\big(\BKcal^\sharp_{j,D}\big)^*$ is the adjoint operator of $\BKcal^\sharp_{j,D} $, that is,
\begin{align*}
\ds \big(\BKcal^\sharp_{j,D}\big)^* [\bphi](x)=\mbox{p.v.}\int_{ \p D}\pd{}{\n(x)} \G_j(x-y) \bphi(y) d\sigma(y) \quad \mbox{a.e. } x \in \p D.
\end{align*}
Moreover,  for  $\bphi \in \mathcal{C}^{1,\alpha}(\p D)$,
\begin{align}\label{conormal-D-sharp}
\pd{\BDcal^\sharp_{j,D}[\bphi]}{\nu_j}\Big|_{+}-\pd{\BDcal^\sharp_{j,D}[\bphi]}{\nu_j}\Big|_{-}=\frac{\p }{\p \ta}\Big(( \bphi\cdot \ta) \n +\frac{\lambda_j}{2\mu_j+\lambda_j}( \bphi\cdot \n ) \ta\Big)\q \mbox{ on }\p D.
\end{align}
\end{lem}
Note that we  drop the $p.v.$ in this stage; this is because $\p D$ is $\mathcal{C}^{2,\eta}$.

Denote by
\begin{align*}
 \ds \mathcal{X}(\p D):=L^2(\p D)^2,\q\mbox{ } \mathcal{X}_{\Psi}(\p D):= L^2(\p D)\times L_{\Psi}^2(\p D),\q\mbox{ } \mathcal{Y}(\p D):= W_{1}^2(\p D)\times L^2(\p D),
\end{align*}
where $W_1^2(\p D)$ is the first $L^2$-Sobolev of space of order $1$ on $\p D$.

The following theorem is of particular importance to us for establishing our representation
formula.
\begin{thm} Suppose that $(\lambda_0-\lambda_2)(\mu_0-\mu_2) \geq 0$  and $0<\lambda_2,\mu_2<\infty$. For any given $(\f_1,\f_2, \g_1, \g_2)\in \mathcal{Y}(\p D) \times \mathcal{Y}(\p D_\ep) $,
there exists a unique solution $(\bvarphi_1,\bvarphi_2,\widetilde{\bpsi}_2,\widetilde{\bvarphi}_0)\in \mathcal{X}(\p D) \times \mathcal{X}(\p D_\ep) $ to the following  integral equations
\begin{align}
  \ds\BScal_{1,D}[\bvarphi_1]\big|_{-}-\BScal_{2,D}[\bvarphi_2]\big|_{+}-\BScal_{2,D_\ep}[\widetilde{\bpsi}_2]&= \f_1 \q\mbox{ on } \p D, \label{eq1}\\
 \nm \ds \pd{\BScal_{1,D}[\bvarphi_1]}{\nu_1}\Big|_{-}-\pd{\BScal_{2,D}[\bvarphi_2]}{\nu_2}\Big|_{+}- \pd{\BScal_{2,D_\ep}[\widetilde{\bpsi}_2]}{\nu_2} &=\f_2 \q\mbox{ on } \p D,\label{eq2} \\
 \nm \ds \BScal_{2,D}[\bvarphi_2]+\BScal_{2,D_\ep}[\widetilde{\bpsi}_2]\big|_{-}- \BScal_{0,D_\ep}[\widetilde{\bvarphi}_0] \big|_{+}&=\g_1 \q\mbox{ on } \p D_\ep, \label{eq3} \\
 \nm \ds\pd{\BScal_{2,D}[\bvarphi_2]}{\nu_2}+ \pd{\BScal_{2,D_\ep}[\widetilde{\bpsi}_2]}{\nu_2}\Big|_{-}-\pd{\BScal_{0,D_\ep}[\widetilde{\bvarphi}_0]}{\nu_0}\Big|_{+}&=\g_2 \q\mbox{ on } \p D_\ep . \label{eq4}
\end{align}
Moreover, if $(\f_2,\g_2)\in L_{\Psi}^2(\p D)\times L_{\Psi}^2(\p D_\ep)$, then $(\bvarphi_2,\widetilde{\bvarphi}_0)\in L_{\Psi}^2(\p D)\times L_{\Psi}^2(\p D_\ep)$.
\end{thm}
\proof The unique solvability of the system of integral equations \eqref{eq1}-\eqref{eq4} can be done  in exactly the same manner  as in \cite[Theorem 5.1]{S}
by using \cite{ES}, see also \cite{book}. By using \eqref{L-Psi}, $\p \BScal_{1,D}[\bvarphi_1]/ \p \nu_1\big|_{-}, \p \BScal_{2,D_\ep}[\widetilde{\bpsi}_2]/\p \nu_2  \in L_{\Psi}^2(\p D)$. It then follows from \eqref{eq2}  that $\p \BScal_{2,D}[\bvarphi_2]/\p \nu_2 \big|_{+}\in L_{\Psi}^2(\p D)$. Since
\begin{align*}
\bvarphi_2=\pd{ \BScal_{2,D}[\bvarphi_2]}{ \nu_2} \Big|_{+}-\pd{ \BScal_{2,D}[\bvarphi_2]}{ \nu_2} \Big|_{-},
\end{align*}
with  $\p \BScal_{2,D}[\bvarphi_2]/\p \nu_2 \big|_{-}\in L_{\Psi}^2(\p D)$, we conclude that $\bvarphi_2 \in L_{\Psi}^2(\p D)$. For any $\eta \in \Psi$, we have
\begin{align*}
\ds \int_{\p D_\ep} \pd{\BScal_{2,D}[\bvarphi_2]}{\nu_2}\cdot \eta d\sigma= \ds \int_{\p D} \pd{\BScal_{2,D}[\bvarphi_2]}{\nu_2}\Big|_{+}\cdot \eta d\sigma=0.
\end{align*}
In order to justify the last equality, we use $\p \BScal_{2,D}[\bvarphi_2]/\p \nu_2 \big|_{+}\in L_{\Psi}^2(\p D)$. Then  $\p \BScal_{2,D}[\bvarphi_2]/\p \nu_2 \in L_{\Psi}^2(\p D_\ep)$. It then follows from \eqref{eq4} that $\p \BScal_{0,D_\ep}[\widetilde{\bvarphi}_0]/ \p \nu_0\big|_{+}\in L_{\Psi}^2(\p D_\ep)$. Thus
\begin{align*}
\ds \widetilde{\bvarphi}_0=\pd{\BScal_{0,D_\ep}[\widetilde{\bvarphi}_0]}{\nu_0}\Big|_{+}-\pd{\BScal_{0,D_\ep}[\widetilde{\bvarphi}_0]}{\nu_0}\Big|_{-}\in L_{\Psi}^2(\p D_\ep).
\end{align*}
This completes the proof.

We now prove a representation theorem for the solution of the transmission
problem \eqref{equation-u} which will be the main ingredient in deriving the  asymptotic
expansion in Theorem \ref{Main-theorem}.
\begin{thm}
The solution $\u_\ep$ to the problem \eqref{equation-u} is  represented by
\begin{equation}\label{representation-u-ep}
\u_\ep(x)=
\left\{
  \begin{array}{ll}
   \ds  \H(x)+\BScal_{0,D_\ep}[\widetilde{\bvarphi}_0](x), & x\in \RR^2\backslash \overline{D}_\ep,  \\
    \nm \ds \BScal_{2,D}[\bvarphi_2](x)+\BScal_{2,D_\ep}[\widetilde{\bpsi}_2](x), & x\in {D}_\ep\backslash \overline{D},  \\
   \nm \ds  \BScal_{1,D}[\bvarphi_1](x), & x\in D,
  \end{array}
\right.
\end{equation}
where  $(\bvarphi_1,\bvarphi_2,\widetilde{\bpsi}_2,\widetilde{\bvarphi}_0)\in \mathcal{X}_{\Psi}(\p D) \times \mathcal{X}_{\Psi}(\p D_\ep) $ is the unique solution to the following integral equations:
\begin{align}
  \ds\BScal_{1,D}[\bvarphi_1]\big|_{-}-\BScal_{2,D}[\bvarphi_2]\big|_{+}-\BScal_{2,D_\ep}[\widetilde{\bpsi}_2]&= 0 \q  \mbox{ on } \p D, \label{eq10}\\
 \nm \ds \pd{\BScal_{1,D}[\bvarphi_1]}{\nu_1}\Big|_{-}-\pd{\BScal_{2,D}[\bvarphi_2]}{\nu_2}\Big|_{+}- \pd{\BScal_{2,D_\ep}[\widetilde{\bpsi}_2]}{\nu_2} &=0 \q \mbox{ on }\p D,\label{eq20} \\
 \nm \ds \BScal_{2,D}[\bvarphi_2]+\BScal_{2,D_\ep}[\widetilde{\bpsi}_2]\big|_{-}- \BScal_{0,D_\ep}[\widetilde{\bvarphi}_0]\big|_{+} &=\H \q  \mbox{ on } \p D_\ep, \label{eq30} \\
 \nm \ds\pd{\BScal_{2,D}[\bvarphi_2]}{\nu_2}+ \pd{\BScal_{2,D_\ep}[\widetilde{\bpsi}_2]}{\nu_2}\Big|_{-}-\pd{\BScal_{0,D_\ep}[\widetilde{\bvarphi}_0]}{\nu_0}\Big|_{+}&=\pd{\H}{\nu_0} \, \mbox{  on } \p D_\ep. \label{eq40}
\end{align}
\end{thm}
\proof Let $(\bvarphi_1,\bvarphi_2,\widetilde{\bpsi}_2,\widetilde{\bvarphi}_0)$ be the unique solution of \eqref{eq10}-\eqref{eq40}. Then clearly $\u_\ep$ defined by   \eqref{representation-u-ep} satisfies the transmission conditions (the conditions on the  fourth and  fifth lines
in \eqref{equation-u}). It is known that
\begin{align*}
\ds \BScal_{0,D_\ep}[\widetilde{\bvarphi}_0](x)=\Gamma_{0}(x)\int_{\p D_\ep}\widetilde{\bvarphi}_0 d\sigma+O(|x|^{-1})\q \mbox{ as }|x|\rightarrow \infty.
\end{align*}
Since $\widetilde{\bvarphi}_0 \in L_{\Psi}^2(\p D_\ep)$, which gives $\int_{\p D_\ep}\widetilde{\bvarphi}_0 d\sigma=0$, and then  $\u_\ep$ defined by   \eqref{representation-u-ep} satisfies $\u_\ep(x)-\H(x)=O(|x|^{-1})$  as $|x|\rightarrow \infty.$ This finishes the proof of the theorem.

We now introduce some notation. Let $\Phi_\ep$ be the diffeomorphism from $\partial D$ onto
$\partial D_\epsilon$ given by $\xe=\Phi_\ep(x) =  x +\ep h(x)\n(x),$
where $x = X(t)\in \p D$. Define the operators $\mathcal{Q}_\ep$ and $\mathcal{R}_\ep$  from
$L^2(\p D) \times L^2(\p D_\ep)$ into $\mathcal{Y}(\p D)$ by
\begin{align}
\ds\mathcal{Q}_\ep(\bvarphi, \widetilde{\bpsi})&:=\bigg(\BScal_{1,D}[\bvarphi]\big|_{-}-\BScal_{0,D_\ep}[\widetilde{\bpsi}]\circ \Phi_\ep\big|_{+},\pd{\BScal_{1,D}[\bvarphi]}{\nu_1}\Big|_{-}- \pd{\BScal_{0,D_\ep}[\widetilde{\bpsi}]}{\nu_0}\circ \Phi_\ep\Big|_{+}\bigg),\label{Q-ep}\\
\nm \ds\mathcal{R}_\ep(\bvarphi, \widetilde{\bpsi})&:=\bigg(\BScal_{2,D}[\bvarphi]\circ \Phi_\ep-\BScal_{2,D}[\bvarphi]\big|_{+}+ \BScal_{2,D_\ep}[\widetilde{\bpsi}]\circ \Phi_\ep\big|_{-}-\BScal_{2,D_\ep}[\widetilde{\bpsi}],\nonumber\\
\nm \ds & \q\q\q\q \q\q \pd{\BScal_{2,D}[\bvarphi]}{\nu_2}\circ \Phi_\ep-\pd{\BScal_{2,D}[\bvarphi]}{\nu_2}\Big|_{+}+ \pd{\BScal_{2,D_\ep}[\widetilde{\bpsi}]}{\nu_2}\circ \Phi_\ep\Big|_{-}-\pd{\BScal_{2,D_\ep}[\widetilde{\bpsi}]}{\nu_2}\bigg),\label{R-ep}
\end{align}
and the matrix-valued function  $\mathcal{H}_\ep$ on $\p D$ by
\begin{align}\label{H-ep}
\ds \mathcal{H}_\ep(x)&:=\bigg(\H\big(x+\ep h(x)\n(x)\big), \pd{\H}{\nu_0}\big(x+\ep h(x)\n(x)\big)\bigg).
\end{align}

By adding \eqref{eq10} to \eqref{eq30} and \eqref{eq20} to \eqref{eq40}, we get  the following lemma.
\begin{lem}\label{lem} \label{representation-formula} Let $(\bvarphi_1,\bvarphi_2,\widetilde{\bpsi}_2,\widetilde{\bvarphi}_0)\in \mathcal{X}_{\Psi}(\p D) \times \mathcal{X}_{\Psi}(\p D_\ep)$  be  the unique solution of \eqref{eq10}-\eqref{eq40}; then  $(\bvarphi_1, \widetilde{\bvarphi}_0)$ and $(\bvarphi_2, \widetilde{\bpsi}_2)$ satisfy  the following system of integral equations:
\begin{align}\label{first-equation}
\ds\mathcal{Q}_\ep(\bvarphi_1, \widetilde{\bvarphi}_0)&=\mathcal{H}_\ep-\mathcal{R}_\ep(\bvarphi_2, \widetilde{\bpsi}_2)\q \mbox{ on }\p D.
\end{align}
\end{lem}

In the next section, we will provide asymptotic expansions of the layer potentials , which are appeared in the system of integral equations \eqref{first-equation} with $(\bvarphi_1, \widetilde{\bvarphi}_0)\in L^2(\p D)\times  L^2(\p D_\ep)$ and $(\bvarphi_2,\widetilde{\bpsi}_2)\in  \mathcal{C}^{1,\eta}(\p D)\times  \mathcal{C}^{1,\eta}(\p D_\ep)$. These asymptotic expansions will help us to derive  the asymptotic expansion of the displacement field $\u_\ep$.

\section{Asymptotic expansions of layer potentials}
Let $a, b \in \RR,$ with $a<b$, and let $X(t): [a,b]\to \RR^2$ be
the arclength parametrization of $\p D$, namely, $X$ is an
$\mathcal{C}^2$-function satisfying $|X'(t)|= 1$ for all $t \in
[a,b]$ and
 $$
 \partial D:=\{x=X(t),  t\in [a,b]\}.
 $$
Then the outward unit normal to $D$, $\n(x)$, is given by
$\n(x)=R_{-\frac{\pi}{2}}X'(t)$, where $R_{-\frac{\pi}{2}}$ is the
rotation by $-{\pi}/{2}$,  the tangential vector at $x$, $\ta(x) =
X'(t)$, and $X'(t)\perp X''(t)$. Set the curvature $\kappa(x)$ to be
defined by
 $$
 X''(t)=\kappa (x)\n(x).
 $$
We will sometimes use $h(t)$ for $h(X(t))$ and $h'(t)$ for the
tangential derivative of $h(x)$.

Then, $\tilde{X}(t)=X(t)+\epsilon h(t)\n(x)=X(t)+\epsilon
h(t)R_{-\frac{\pi}{2}}X'(t)$ is a parametrization of $\p {D}_\ep$. By
${\n} (\tilde{x})$, we denote the outward unit normal to
$\p D_\ep$ at $\tilde{x}$.

Let $\xe=x+\ep h(x)\n(x)\in \p D_\ep$ for $x\in \p D$. The following asymptotic  expansions of $\n(\xe)$ and the length element
$d{\sigma_\epsilon}(\tilde{x})$ hold \cite{AKLZ1}:
 \begin{align} \label{nuexpan}
 \ds \n(\tilde{x}) = \n(x) - \ep h'(t) \ta (x) + O(\ep^2),
 \end{align}
 and
 \begin{align} \label{sigexp}
\ds  d{\sigma_\epsilon}(\tilde{x}) = \big(1 -\ep \kappa(x)h(x)+ O(\ep^2)\big)
 d\sigma(x).
 \end{align}
Here, the remainder term $O(\ep^2)$ is bounded by $C\ep^2$ for some constant $C$ which depends only  on $\mathcal{C}^2$-norm of $\p D$
and $\mathcal{C}^1$-norm of $h$.

Let $\bphi (x)$ and $\phi(x)$ be a vector function and scalar function, respectively,  which belong to $ \mathcal
{C}^2([a,b])$ for $x=X(\cdot)\in \p D$. By ${d}/{dt}$,  we denote the tangential derivative in the
direction of $\ta(x)=X'(t)$.  We have
\begin{align*}
\ds \frac{d}{dt}\big(\bphi(x)\big)=\nabla \bphi(x) X'(t)=\frac{\p \bphi}{\p
\ta}(x),\q\q\frac{d}{dt}\big(\phi(x)\big)=\nabla \phi(x) \cdot X'(t)=\frac{\p \phi}{\p
\ta}(x).
\end{align*}

The  restriction of the Lam\'e system $\Lcal_{\lambda_j,\mu_j}$ in $D$ to
a neighborhood  of $\p D$  can be expressed as follows \cite{JFH}:
\begin{align}\label{local-Lame}
\ds \Lcal_{\lambda_j,\mu_j}\bphi(x)=& \mu_j \pdl{\bphi}{\n^2}(x)+\lambda_j \nabla\nabla\cdot \bphi(x)\cdot \n(x)\n(x)+\mu_j \nabla(\nabla \bphi)^T(x)\n(x)\n(x)\nonumber\\
\nm\ds&-\kappa(x)\pd{ \bphi}{\nu_j}(x)+\frac{d}{dt}\Big(\big(\mathbb{C}_j\widehat{\nabla} \bphi(x) \big)\ta(x)\Big),\q x\in \p D.
\end{align}

We have from \cite{JFH} the following lemma.
\begin{lem}
Let $\widetilde{\bphi}  \in L^2(\p D_\ep)$, we denote  by $\bphi:= \widetilde{\bphi}\circ \Phi_\ep$. For $j=0,2$, the following  asymptotic expansions hold:
\begin{align}
\ds\BScal_{j,D_\ep}[\widetilde{\bphi}]\circ \Phi_\ep\big|_{\pm}&= \BScal_{j,D}[\bphi]-\ep \BScal_{j,D}[\kappa h \bphi]
+\ep \Big(h\pd{\BScal_{j,D}[\bphi]}{\n}+\BDcal^{\sharp}_{j ,D}[h\bphi]\Big)\Big|_{\pm}\nonumber\\
\nm \ds &\q + O_1(\ep^2) \q \mbox{on }\p D,\label{asym-3}\\
\nm \ds \pd{\BScal_{j,D_\ep}[\widetilde{\bphi}]}{\nu_j}\circ \Phi_\ep\Big|_{\pm}&=\pd{\BScal_{j,D}[\bphi]}{\nu_j}\Big|_{\pm}+\ep \bigg (\kappa h\pd{\BScal_{j,D}[\bphi]}{\nu_j}-\pd{\BScal_{j,D}[\kappa h\bphi]}{\nu_j}\bigg)\bigg|_{\pm}\nonumber \\
\nm \ds &\q+\ep\bigg(\pd{\BDcal^{\sharp}_{j,D}[h\bphi]}{\nu_j}- \frac{\p}{\p \ta}\Big(h\big(\mathbb{C}_j\widehat{\nabla }\BScal_{j,D}[\bphi]\big)\ta\Big)\bigg)\bigg|_{\pm}+O_2(\ep^2)\q \mbox{on }\p D,\label{asym-1}
\end{align}
where  $\|O_1(\ep^2)\|_{W_{1}^2(\p D)}, \|O_2(\ep^2)\|_{L^2(\p D)}\leq C\ep^2$ for some constant $C$ depends only on $\lambda_j,\mu_j$, the $\mathcal{C}^2$-norm of $X$, and the $\mathcal{C}^1$-norm of $h$.
\end{lem}

Now, we are going to derive  asymptotic expansions of $\BScal_{2,D}[\bphi](\xe)$ and $ {\p \BScal_{2,D}[\bphi]}/{\p \nu_2}(\xe)$  for $\bphi \in \mathcal{C}^{1,\eta}(\p D)$, $\xe=x+\ep h(x)\n(x) \in \p D_\ep$, and $x\in \p D$. Because $\p D$ is $C^{2, \eta}$, $\BScal_{2,D}[\bphi]$  is $\mathcal{C}^{2,\eta}(\RR^2\backslash D)$. We have
\begin{align}
\ds&\Big| \nabla \BScal_{2,D}[\bphi](\xe)- \nabla \BScal_{2,D}[\bphi](x)\big|_{+}-\ep h(x) \nabla^2 \BScal_{2,D}[\bphi](x)\n(x)\big|_{+}\Big|\nonumber\\
\nm \ds &\q =\Big|\int_{0}^1 \Big[\nabla^2 \BScal_{2,D}[\bphi]\big(x+l(\xe-x)\big)-\nabla^2 \BScal_{2,D}[\bphi](x)\big|_{+}\Big]dl(\xe-x)\Big|\nonumber\\
\nm \ds &\q  \leq C\big|\xe -x\big|^{1+\eta}\big\|\BScal_{2,D}[\bphi]\big\|_{\mathcal{C}^{2,\eta}(\RR^2\backslash D)}\nonumber \\
\nm \ds &\q   \leq C  \ep^{1+\eta}\big\|\bphi\big\|_{\mathcal{C}^{1,\eta}(\p  D)}.\label{holder-estimate}
\end{align}
Thus
\begin{align}\label{eqq1}
 \ds \pd{\BScal_{2,D}[\bphi]}{\ta}(\xe)& =\Big(\nabla \BScal_{2,D}[\bphi](x)\big|_{+}+\ep h(x) \nabla^2 \BScal_{2,D}[\bphi](x)\n(x)\big|_{+}+O(\ep^{1+\eta})\Big)\nonumber\\
 \nm \ds &\q\q\q\q\q\q\q \q\q\q \q\q\q \times \Big(\ta(x)+\ep h'(t)\n(x)+O(\ep^2)\Big)\nonumber\\
 \nm \ds &=\frac{\p }{\p \ta}\Big(\BScal_{2,D}[\bphi](x)+\ep h(x)\frac{\p\BScal_{2,D}[\bphi]}{ \p \n}\Big|_{+}(x)\Big)+O(\ep^{1+\eta}), \q x\in \p D,
\end{align}
where $\|O(\ep^{1+\eta})\|_{L^2(\p D)}$ is bounded by $C \ep^{1+\eta} \|\bphi\|_{\mathcal{C}^{1,\eta}(\p D)}$. Similarly to \eqref{holder-estimate}, we get
\begin{align}\label{eqq2}
\ds  \BScal_{2,D}[\bphi](\xe)&= \BScal_{2,D}[\bphi](x)+\ep h(x)\frac{\p\BScal_{2,D}[\bphi]}{ \p \n}\Big|_{+}(x)+O(\ep^{2}),\q x\in \p D,
\end{align}
where $\|O(\ep^{2})\|_{L^2(\p D)}$ is bounded by $C \ep^{2} \|\bphi\|_{\mathcal{C}^{1,\eta}(\p D)}$. In conclusion, we get  from \eqref{eqq1} and \eqref{eqq2}
that
\begin{align}
\ds  \BScal_{2,D}[\bphi](\xe)&= \BScal_{2,D}[\bphi](x)+\ep h(x)\pd{\BScal_{2,D}[\bphi]}{\n}\Big|_{+}(x)+O(\ep^{1+\eta}), \q x\in \p D,\label{asym-6}
\end{align}
where $\|O(\ep^{1+\eta})\|_{W_1^2(\p D)}$ is bounded by $C \ep^{1+\eta} \|\bphi\|_{\mathcal{C}^{1,\eta}(\p D)}$.

It follows from \eqref{conormal-derivative},  \eqref{nuexpan}, \eqref{local-Lame},   and Taylor expansion that, for $\xe\in  \p D$,
\begin{align}
\ds \pd{\BScal_{2,D}[\bphi]}{\nu_2}(\xe)=&\lambda_2 \nabla \cdot \BScal_{2,D}[\bphi](\xe)\n(\xe)+\mu_2\Big( \nabla \BScal_{2,D}[\bphi](\xe)+\nabla \BScal^T_{2,D}[\bphi](\xe)\Big)\n(\xe)\nonumber\\
\nm \ds= &\Big[\lambda_2 \nabla \cdot \BScal_{2,D}[\bphi](x)+\mu_2\Big( \nabla \BScal_{2,D}[\bphi](x)+\nabla \BScal^T_{2,D}[\bphi](x)\Big)\Big]\Big|_{+}\n(x)\nonumber\\
\nm \ds &-\ep h'(t)\Big[\lambda_2 \nabla \cdot \BScal_{2,D}[\bphi](x)+\mu_2\Big( \nabla \BScal_{2,D}[\bphi](x)+\nabla \BScal^T_{2,D}[\bphi](x)\Big)\Big]\Big|_{+}\ta(x)\nonumber\\
\nm \ds &+\ep h(x)\Big[\lambda_2\nabla \nabla \cdot \BScal_{2,D}[\bphi](x)\cdot \n(x)\n(x)+\mu_2 \nabla \nabla \BScal_{2,D}[\bphi](x)\n(x)\n(x)\nonumber\\
\nm \ds &\q\q\q\q  +\mu_2\nabla \nabla \BScal^T_{2,D}[\bphi](x)\n(x)\n(x)\Big]\Big|_{+}+O(\ep^{1+\eta})\nonumber\\
\nm\ds&=\pd{\BScal_{2,D}[\bphi]}{\nu_2}\Big|_{+}(x)+\ep \kappa(x)h(x)\pd{\BScal_{2,D}[\bphi]}{\nu_2}\Big|_{+}(x)\nonumber\\
\nm \ds &\q -\ep \frac{d}{dt}\bigg(h(x)\Big(\mathbb{C}_2\widehat{\nabla }\BScal_{2,D}[\bphi](x)\Big)\ta(x)\bigg)\bigg|_{+}+O(\ep^{1+\eta}),\q x\in \p D,\label{asym-5}
\end{align}
where $\|O(\ep^{1+\eta})\|_{L^2(\p D)}$ is bounded by $C \ep^{1+\eta} \|\bphi\|_{\mathcal{C}^{1,\eta}(\p D)}$. In order to justify the last equality in \eqref{asym-5}, we use the representation of the Lam\'e system on $\p D$ given in \eqref{local-Lame}.

We now expand $\BScal_{2,D_\ep}[\widetilde{\bphi}](x)$ and $\p \BScal_{2,D_\ep}[\widetilde{\bphi}]/\p\nu_2(x)$ for $x\in \p D$ when $\widetilde{\bphi} \in \mathcal{C}^{1,\eta}(\p D_\ep)$. Let $\f$ be a $\mathcal{C}^{1,\eta}$  vector function on $\p D$ and let  $\w$ be the solution to $\Lcal_{\lambda_2, \mu_2}\w=0$  in $D$ satisfying $\w=\f$ on $\p D$. Then, we get
\begin{align}\label{eq1000}
\ds  \int_{\p D}\pd{\BScal_{2,D_\ep}[\widetilde{\bphi}]}{\nu_2}(x)\cdot \f(x) d\sigma(x)&= \int_{\p D}\BScal_{2,D_\ep}[\widetilde{\bphi}](x)\cdot \pd{\w}{\nu_2}(x) d\sigma(x)\nonumber\\
\nm \ds &= \int_{\p D_\ep}\widetilde{\bphi}(\xe)\cdot \BScal_{2,D}\Big[\pd{\w}{\nu_2}\Big](\xe) d\sigma(\xe).
\end{align}
Define $\bphi:=\widetilde{\bphi}\circ \Phi_\ep$. By using  \eqref{sigexp}, we get
\begin{align*}
\ds \int_{\p D}\pd{\BScal_{2,D_\ep}[\widetilde{\bphi}]}{\nu_2}\cdot\f d\sigma&= \int_{\p D}\bphi\cdot\bigg(\BScal_{2,D}\Big[\pd{\w}{\nu_2}\Big]+ \ep h\pd{\BScal_{2,D}}{\n}\Big[\pd{\w}{\nu_2}\Big]\Big |_{+}+O(\ep^{1+\eta}) \bigg)\\
\nm\ds &\q \q \q\q \q \q \q\q \q \q\q\q \q \q\q\q \q  \times \Big(1-\ep \kappa h+O(\ep^2)\Big)d\sigma\\
\nm\ds &   =\int_{\p D}\Big(\BScal_{2,D}[\bphi]+\ep \BDcal^{\sharp}_{2,D}[h\bphi]\big|_{-} -\ep \BScal_{2,D}[\kappa h \bphi] \Big)\cdot\pd{\w}{\nu_2}d\sigma+O(\ep^{1+\eta})\\
\nm\ds &=\int_{\p D}\bigg(\pd{\BScal_{2,D}[\bphi]}{\nu_2}\Big|_{-}+\ep \pd{\BDcal^{\sharp}_{2,D}[h\bphi]}{\nu_2}\Big|_{-} -\ep
\pd{\BScal_{2,D}[\kappa h \bphi]}{\nu_2}\Big|_{-} \bigg)\cdot\f d\sigma\\
\nm \ds &\q+O(\ep^{1+\eta}).
\end{align*}
Therefore, the following asymptotic expansion holds:
\begin{align}
\ds \pd{\BScal_{2,D_\ep}[\widetilde{\bphi}]}{\nu_2}=\pd{\BScal_{2,D}[\bphi]}{\nu_2}\Big|_{-}+\ep \pd{\BDcal^{\sharp}_{2,D}[h\bphi]}{\nu_2}\Big|_{-} -\ep
\pd{\BScal_{2,D}[\kappa h \bphi]}{\nu_2}\Big|_{-}+O(\ep^{1+\eta})\q \mbox{on }\p D,\label{asym-7}
\end{align}
where the  remainder term $O(\ep^{1+\eta})$ is in $L^2(\p D)$.

Let $\bphi=\widetilde{\bphi}\circ \Phi_\ep$ for $\widetilde{\bphi}\in \mathcal{C}^{1,\eta}(\p D_\ep)$. Let $\f$ be a $\mathcal{C}^{1,\eta}$ vector function on $\p D$. Similarly to \eqref{eq1000},  we have
\begin{align*}
\ds \int_{\p D}\pd{\BScal_{2,D_\ep}[\widetilde{\bphi}]}{\ta}\cdot\f d\sigma&= -\int_{\p D}\bphi\cdot\bigg(\BScal_{2,D}\Big[\pd{\f}{\ta}\Big]+ \ep h\pd{\BScal_{2,D}}{\n}\Big[\pd{\f}{\ta}\Big]\Big |_{+}+O(\ep^{1+\eta}) \bigg)\\
\nm\ds &\q \q \q\q \q \q \q\q \q \q\q\q \q \q\q\q\q\q \times \Big(1-\ep \kappa h+O(\ep^2)\Big)d\sigma \\
\nm\ds &   =-\int_{\p D}\Big(\BScal_{2,D}[\bphi]+\ep \BDcal^{\sharp}_{2,D}[h\bphi]\big|_{-} -\ep \BScal_{2,D}[\kappa h \bphi] \Big)\cdot\pd{\f}{\ta}d\sigma+O(\ep^{1+\eta})\\
\nm\ds &=\int_{\p D}\bigg(\pd{\BScal_{2,D}[\bphi]}{\ta}+\ep \pd{\BDcal^{\sharp}_{2,D}[h\bphi]}{\ta}\Big|_{-} -\ep
\pd{\BScal_{2,D}[\kappa h \bphi]}{\ta} \bigg)\cdot\f d\sigma+O(\ep^{1+\eta}).
\end{align*}
Thus
\begin{align}\label{eqq3}
\ds &\pd{\BScal_{2,D_\ep}[\widetilde{\bphi}]}{\ta}=\pd{\BScal_{2,D}[\bphi]}{\ta}+\ep \pd{\BDcal^{\sharp}_{2,D}[h\bphi]}{\ta}\Big|_{-}-\ep
\pd{\BScal_{2,D}[\kappa h \bphi]}{\ta}+O(\ep^{1+\eta}) \q \mbox{ on } \p D,
\end{align}
where the  remainder term $O(\ep^{1+\eta})$ is in $L^2(\p D)$. In a similar way, we get
\begin{align}\label{eqq4}
\ds\BScal_{2,D_\ep}[\widetilde{\bphi}]&= \BScal_{2,D}[\bphi]+\ep \BDcal^{\sharp}_{2,D}[h\bphi]\big |_{-}-\ep \BScal_{2,D}[\kappa h\bphi]+O(\ep^{2})\q
 \mbox{ on } \p D,
\end{align}
with  the  remainder term $O(\ep^{2})$ is in $L^2(\p D)$. It then follows from \eqref{eqq3} and \eqref{eqq4} that
\begin{align}
\ds\BScal_{2,D_\ep}[\widetilde{\bphi}]&= \BScal_{2,D}[\bphi]+\ep \BDcal^{\sharp}_{2,D}[h\bphi]\big |_{-}-\ep \BScal_{2,D}[\kappa h\bphi]+O(\ep^{1+\eta})\q
\mbox{ on } \p D,\label{asym-8}
\end{align}
where  the  remainder term $O(\ep^{1+\eta})$ is in $W_1^2(\p D)$.

The following proposition is a direct consequence of \eqref{Q-ep}, \eqref{R-ep},  \eqref{asym-3}, \eqref{asym-1}, \eqref{asym-6}, \eqref{asym-5}, \eqref{asym-7}, and \eqref{asym-8}.
\begin{prop}\label{prop} The following expansions hold on $\p D$:
\begin{align*}
\ds &\mathcal{Q}_\ep(\bvarphi, \widetilde{\bpsi})=\mathcal{Q}_0(\bvarphi,\bpsi)-\ep \mathcal{Q}_1(\bpsi) +O(\ep^2)\q \mbox{for } (\bvarphi, \widetilde{\bpsi})\in L^2(\p D)\times  L^2(\p D_\ep),\\
\nm \ds &\mathcal{R}_\ep(\bvarphi, \widetilde{\bpsi})=\ep \mathcal{R}_1(\bvarphi,\bpsi) +o(\ep)\q \mbox{ for } (\bvarphi, \widetilde{\bpsi})\in \mathcal{C}^{1,\eta}(\p D)\times  \mathcal{C}^{1,\eta}(\p D_\ep),
\end{align*}
where $\bpsi=\widetilde{\bpsi}\circ \Phi_\ep$, the remainder terms $O(\ep^2)$ and $o(\ep)$ are in $W^2_1(\p D)\times L^2(\p D)$, and the operators
 $\mathcal{Q}_0: \mathcal{X}(\p D)\rightarrow \mathcal{Y}(\p D)$,  $\mathcal{Q}_1: L^2(\p D)\rightarrow \mathcal{Y}(\p D)$, and  $\mathcal{R}_1: \mathcal{C}^{1,\eta}(\p D)\times  \mathcal{C}^{1,\eta}(\p D)\rightarrow \mathcal{Y}(\p D)$ are defined  by
\begin{align}
\ds& \mathcal{Q}_0(\bvarphi,\bpsi)=\bigg(\BScal_{1,D}[\bvarphi]-\BScal_{0,D}[\bpsi],
\pd{\BScal_{1,D}[\bvarphi]}{\nu_1}\Big|_{-}-\pd{\BScal_{0,D}[\bpsi]}{\nu_0}\Big|_{+} \bigg),\label{Q0}\\
\nm \ds &\mathcal{Q}_1(\bpsi)=\bigg(-\BScal_{0,D}[\kappa h \bpsi]+ h \pd{\BScal_{0,D}[\bpsi]}{\n}\Big|_{+}+   \BDcal^{\sharp}_{0,D}[h \bpsi]\Big|_{+},\label{Q1}\nonumber\\
\nm \ds &\q \q \q\,\,\,\,  \kappa h \pd{\BScal_{0,D}[\bpsi]}{\nu_0}\Big|_{+}-\pd{\BScal_{0,D}[\kappa\ h \bpsi]}{\nu_0}\Big|_{+} + \pd{\BDcal^{\sharp}_{0,D}[h\bpsi]}{\nu_0}\Big|_{+}- \frac{\p }{\p \ta}\Big(h\big(\mathbb{C}_0\widehat{\nabla }\BScal_{0,D}[\bpsi]\big)\ta \Big)\Big|_{+}\bigg),\\
\nm \ds &\mathcal{R}_1(\bvarphi,\bpsi)=\bigg(h \pd{\BScal_{2,D}[\bvarphi]}{\n}\Big|_{+}+h \pd{\BScal_{2,D}[\bpsi]}{\n}\Big|_{-},\kappa h \pd{\BScal_{2,D}[\bvarphi]}{\nu_2}\Big|_{+}+\kappa h \pd{\BScal_{2,D}[\bpsi]}{\nu_2}\Big|_{-}\nonumber\\
\nm  \ds &\q \q\q\q \q\q\q\q\q\q\q\q\q- \frac{\p }{\p \ta}\Big(h\big(\mathbb{C}_2\widehat{\nabla }\BScal_{2,D}[\bvarphi]\big)\ta \Big)\Big|_{+}- \frac{\p }{\p \ta}\Big(h\big(\mathbb{C}_2\widehat{\nabla }\BScal_{2,D}[\bpsi]\big)\ta \Big)\Big|_{-} \bigg).\label{R1}
\end{align}
\end{prop}

The following proposition holds.
\begin{prop} Let $(\bvarphi_1,\widetilde{\bvarphi}_0)\in L^2(\p D)\times  L^2(\p D_\ep)$ be the solution of \eqref{eq10}-\eqref{eq40}. Then the following asymptotic expansion holds:
 \begin{align}\label{phi-0-phi-1}
 \ds  \mathcal{Q}_0(\bvarphi_1,\bvarphi_0)-\ep \big[\mathcal{Q}_1(\bvarphi_0)-  \mathcal{Z}(\bvarphi_1)\big]=\mathcal{H}_\ep+ o(\ep) \q \mbox{ on } \p D,
\end{align}
where $\bvarphi_0=\widetilde{\bvarphi}_0\circ \Phi_{\ep}$, the remainder term $o(\ep)$ is in $W^2_1(\p D)\times L^2(\p D)$, $\mathcal{H}_\ep$ is defined by \eqref{H-ep}, $\mathcal{Q}_0$  and  $\mathcal{Q}_1$ are  defined in \eqref{Q0} and \eqref{Q1}, respectively, and the operator $\mathcal{Z}$ is  defined from $L^2(\p D)$ into $\mathcal{Y}(\p D)$ by
\begin{align}\label{Z}
\ds  \mathcal{Z}(\bvarphi_1)&:=\bigg(h \pd{\BScal_{1, D}[\bvarphi_1]}{\n}\Big|_{-}+h\Big(\mathbb{K}_{2,1}\widehat{\nabla}\BScal_{1, D}[\bvarphi_1]\Big)\n\Big|_{-},\nonumber \\
\nm\ds &\q\q\q\q\q\q\q\q \q\q\q\kappa h\pd{\BScal_{1, D}[\bvarphi_1]}{\nu_1}\Big|_{-}- \frac{\p}{\p \ta }\Big(h\big(\mathbb{M}_{2,1}\widehat{\nabla }\BScal_{1, D}[\bvarphi_1] \big)\ta \Big)\Big|_{-}\bigg),
\end{align}
with
\begin{align*}
\ds \mathbb{M}_{2,1}&=\frac{\lambda_2 (\lambda_1+2\mu_1)}{\lambda_2+2\mu_2} \I \otimes \I+2 \mu_1 \mathbb{I}+ \frac{4(\mu_2-\mu_1)(\lambda_2+\mu_2)}{\lambda_2+2\mu_2} \I\otimes (\ta \otimes \ta ),\\
\nm \ds \mathbb{K}_{2,1}&=\frac{(\lambda_1-\lambda_2)\mu_2+2(\mu_2-\mu_1) (\lambda_2+\mu_2)}{\mu_2(\lambda_2+2\mu_2)} \I \otimes \I+2\big(\frac{\mu_1}{\mu_2}-1\big) \mathbb{I}\\
\nm\ds &\q + \frac{2(\mu_1-\mu_2)(\lambda_2+\mu_2)}{\mu_2(\lambda_2+2\mu_2)} \I\otimes (\ta \otimes \ta).
\end{align*}
\end{prop}
\proof
It follows from  the lemma  \ref{lem} and the  proposition \ref{prop} that
\begin{align}
\ds \mathcal{Q}_0(\bvarphi_1,\bvarphi_0)- \ep\big[\mathcal{Q}_0(\bvarphi_0)-  \mathcal{R}_1(\bvarphi_2, \bpsi_2)\big]=\mathcal{H}_\ep+o(\ep)\q \mbox{ on }\p D, \label{eq100}
\end{align}
where  the remainder term $o(\ep)$ is in $W^2_1(\p D)\times L^2(\p D)$ and
\begin{align*}
 \ds\mathcal{R}_1(\bvarphi_2, \bpsi_2)&=\bigg( h \pd{\BScal_{2,D}[\bvarphi_2]}{\n}\Big|_{+}+h \pd{\BScal_{2,D}[\bpsi_2]}{\n}\Big|_{-},
\kappa h \pd{\BScal_{2,D}[\bvarphi_2]}{\nu_2}\Big|_{+}+\kappa h \pd{\BScal_{2,D}[\bpsi_2]}{\nu_2}\Big|_{-}\\
\nm\ds &\q\q\q\q\q\q\q
-\frac{\p}{\p \ta}\Big(h \big(\mathbb{C}_2 \widehat{\nabla} \BScal_{2,D}[\bvarphi_2]\big)\ta\Big)\Big|_{+}-\frac{\p}{\p \ta}\Big(h \big(\mathbb{C}_2 \widehat{\nabla} \BScal_{2,D}[\bpsi_2]\big)\ta\Big)\Big|_{-}\bigg).
\end{align*}
Let  $(\bvarphi_1, \bvarphi_2, \widetilde{\bpsi}_2)\in L^2(\p D)\times \mathcal{C}^{1,\eta}(\p D)\times \mathcal{C}^{1,\eta}(\p D_\ep)$ be the solution of \eqref{eq10}-\eqref{eq40}. According to Appendix, $(\bvarphi_1, \bvarphi_2, \widetilde{\bpsi}_2)$ satisfies the following equations along the interface $\p D$:
\begin{align*}
\left\{
  \begin{array}{lll}
   \ds \pd{\BScal_{2,D}[\bvarphi_2]}{\nu_2}\Big |_{+}+\pd{\BScal_{2,D_\ep}[\widetilde{\bpsi}_2]}{\nu_2}=\pd{\BScal_{1,D}[\bvarphi_1]}{\nu_1}\Big |_{-},\\
    \nm \ds  \nabla \BScal_{2,D}[\bvarphi_2]\,\n\big |_{+}+\nabla \BScal_{2,D_\ep}[\widetilde{\bpsi}_2]\,\n=\nabla \BScal_{1,D}[\bvarphi_1]\, \n\big |_{-}+\big(\mathbb{K}_{2,1}\nabla \BScal_{1,D}[\bvarphi_1]\big)\n\big |_{-},\\
\nm \ds \Big(\mathbb{C}_2 \widehat{\nabla} \BScal_{2,D}[\bvarphi_2]\Big)\ta\Big |_{+}+\Big(\mathbb{C}_2 \widehat{\nabla} \BScal_{2,D_\ep}[\widetilde{\bpsi}_2]\Big)\ta=\Big(\mathbb{M}_{2,1} \widehat{\nabla} \BScal_{1,D}[\bvarphi_1]\Big)\Big |_{-}\ta.
  \end{array}
\right.
\end{align*}
Let $\bpsi_2=\widetilde{\bpsi}_2\circ \Phi_\ep$. The following asymptotic expansions hold on $\p D$:
\begin{align}
\left\{
  \begin{array}{lll}
   \ds \pd{\BScal_{2,D}[\bvarphi_2]}{\nu_2}\Big |_{+}+\pd{\BScal_{2,D}[\bpsi_2]}{\nu_2}\Big |_{-}+o(1)=\pd{\BScal_{1,D}[\bvarphi_1]}{\nu_1}\Big |_{-},\\
    \nm \ds  \nabla \BScal_{2,D}[\bvarphi_2]\,\n\big |_{+}+\nabla \BScal_{2,D}[\bpsi_2]\,\n\big |_{-}+o(1)=\nabla \BScal_{1,D}[\bvarphi_1]\, \n\big |_{-}+\big(\mathbb{K}_{2,1}\nabla \BScal_{1,D}[\bvarphi_1]\big)\n\big |_{-},\\
\nm \ds \Big(\mathbb{C}_2 \widehat{\nabla} \BScal_{2,D}[\bvarphi_2]\Big)\ta\Big |_{+}+\Big(\mathbb{C}_2 \widehat{\nabla} \BScal_{2,D}[{\bpsi}_2]\Big)\ta\Big |_{-}+o(1)=\Big(\mathbb{M}_{2,1} \widehat{\nabla} \BScal_{1,D}[\bvarphi_1]\Big)\Big |_{-}\ta,\label{system}
  \end{array}
\right.
\end{align}
where $o(1)$ is in $W_1^2(\p D)$, which implies  $ \mathcal{R}_1(\bvarphi_2, \bpsi_2)=\mathbb{Z}(\bvarphi_1)+o(1)$ on $\p D$. Therefore, \eqref{eq100} holds, as claimed.


\section{Asymptotic expansion of the displacement field}
The following lemma is important to us.
\begin{lem} Suppose that $(\lambda_0-\lambda_1)(\mu_0-\mu_1)\geq 0$ and $0<\lambda_1,\mu_1<\infty.$ For any given $(\F,\GG)\in \mathcal{Y}(\p D)$,
there exists a unique pair $(\f,\g)\in \mathcal{X}(\p D)$ such that
\begin{equation}\label{solvability}
\mathcal{Q}_0(\f,\g)-\ep \big[\mathcal{Q}_1(\f)-  \mathcal{Z}(\g)\big]=\big(\F, \GG\big).
\end{equation}
Furthermore, there exists a constant $C$ depending only on  $\lambda_0, \mu_0, \lambda_1, \mu_1$,
and the Lipschitz character of $D$ such that
\begin{align}\label{estimate}
\|\f\|_{L^2(\p D)}+\|\g\|_{L^2(\p D)}\leq C \Big(\|\F\|_{W_{1}^2(\p D)}+\|\GG\|_{L^2(\p D)}\Big).
\end{align}
\end{lem}
\proof The operator $\mathcal{Q}: \mathcal{X}(\p D) \rightarrow \mathcal{Y}(\p D)$ defined by $\mathcal{Q}(\f,\g)= \mathcal{Q}_1(\f)-  \mathcal{Z}(\g)$ is bounded on $\mathcal{X}(\p D)$.
It is proved in \cite{ES} that the operator $\mathcal{Q}_0: \mathcal{X}(\p D) \rightarrow \mathcal{Y}(\p D)$ is invertible. For $\ep$ small enough, it follows from Theorem $1.16$, section $4$  of \cite{Kato}, that  the operator $\mathcal{Q}_0-\ep\mathcal{Q}$ is invertible. This completes the proof of solvability of \eqref{solvability}. The estimate \eqref{estimate} is a consequence of solvability and the closed graph  theorem.

\subsection{Proof of the theorem \ref{Main-theorem}}\label{proof}
For $\xe=x+\ep h(x)\n(x)\in \p D_\ep$. We have the following Taylor expansion
\begin{align}
\ds \H(\tilde x)=\H(x)+\ep h(x)\pd{\H}{\n}(x)+O(\ep^2), \q x\in \p D.\label{Expansion-H-xe}
\end{align}
Similarly, by the Taylor expansion, \eqref{nuexpan}, and \eqref{local-Lame},  we obtain
\begin{align}
\ds&\pd{\H}{\nu_0}(\xe)=\pd{\H}{\nu_0}(x)+\ep \kappa(x)h(x)\pd{\H}{\nu_0}(x)-\ep \frac{\p }{\p \ta }\Big(h\big(\mathbb{C}_0\widehat{\nabla }\H\big)\ta\Big)(x)+O(\ep^2),\q x\in \p D.\label{Expansion-conormal-H-xe}
\end{align}
It then follows from \eqref{H-ep}, \eqref{Expansion-H-xe},  and \eqref{Expansion-conormal-H-xe} that
\begin{align}\label{expansion-H}
\ds \mathcal{H}_\ep&=\Big(\H , \pd{\H}{\nu_0}\Big)+\ep \Big(h \pd{\H}{\n},\kappa h\pd{\H}{\nu_0}-\frac{\p}{\p \ta}\big(h\big[\mathbb{C}_0\widehat{\nabla} \H\big]\ta \big) \Big)+O(\ep^2)\nonumber\\
\nm\ds &:=\mathcal{H}_0+\ep\mathcal{H}_1+O(\ep^2)\q \mbox{ on } \p D.
\end{align}
Now, we introduce  $(\bvarphi^0_1,\bvarphi^0_0)$ and $(\bvarphi^1_1,\bvarphi^1_0)$  by the following recursive relations
\begin{align}
\ds \mathcal{Q}_0(\bvarphi^0_1,\bvarphi^0_0)&=\mathcal{H}_0,\label{term-0}\\
\nm\ds   \mathcal{Q}_0(\bvarphi^1_1,\bvarphi^1_0)&=\mathcal{H}_1+ \mathcal{Q}_1(\bvarphi^0_0)-  \mathcal{Z}(\bvarphi_1^0)\label{term-1},
\end{align}
 where  $\mathcal{Q}_0$, $\mathcal{Q}_1$, and $\mathcal{Z}$ are  defined in \eqref{Q0},  \eqref{Q1}, and \eqref{Z}, respectively.  One can see the existence and uniqueness of $(\bvarphi^n_1,\bvarphi^n_0)$ for  $n=0,1$, by using \cite{ES}.

 Let $(\bvarphi_1,\widetilde{\bvarphi}_0)$ be the solution of \eqref{eq10}-\eqref{eq40}. It follows from \eqref{term-0} and \eqref{term-1} that
\begin{align}\label{phi-ep-psi-ep}
\ds &\mathcal{Q}_0(\bvarphi_1-\bvarphi^0_1-\ep \bvarphi^1_1,\widetilde{\bvarphi}_0\circ \Phi_\ep-\bvarphi^0_0-\ep \bvarphi^1_0)-\ep \big[\mathcal{Q}_1(\widetilde{\bvarphi}_0\circ \Phi_\ep-\bvarphi^0_0-\ep \bvarphi^1_0)-  \mathcal{Z}(\bvarphi_1-\bvarphi^0_1-\ep \bvarphi^1_1)\big]\nonumber\\
\nm\ds &\q =\mathcal{H}_\ep-\mathcal{H}_0-\ep\mathcal{H}_1+ o(\ep) \q \mbox{ on } \p D,
\end{align}
where $\|o(\ep)\|_{W_1^2(\p D)\times L^2(\p D)}\leq C\ep^{1+\eta}$  for some $\eta>0$ and $(\bvarphi_1^0, \bvarphi_0^0)$ and $(\bvarphi_1^1,\bvarphi_0^1)$  are the solutions to \eqref{term-0} and  \eqref{term-1}, respectively.

The following lemma follows immediately from  \eqref{expansion-H}, \eqref{phi-ep-psi-ep},  and the estimate   in \eqref{estimate}.
\begin{lem}\label{Lemma} Let $(\bvarphi_1,\widetilde{\bvarphi}_0)$ be the solution of \eqref{eq10}-\eqref{eq40}.  For $\ep$ small enough, there exists $C$ depending only on   $(\lambda_j, \mu_j)$ for j=0,1,2, the $\mathcal{C}^2$-norm of $X$, and the $\mathcal{C}^1$-norm of $h$ such that
\begin{equation}\label{estimate-important}
\ds \Big\| \bvarphi_1-\bvarphi^0_1-\ep \bvarphi^1_1\Big\|_{L^2(\p D)}+
\Big\|\widetilde{\bvarphi}_0\circ \Phi_\ep-\bvarphi^0_0-\ep \bvarphi^1_0\Big\|_{L^2(\p D)}\leq C \ep^{1+\eta}
\end{equation}
for some $\eta>0$, where $(\bvarphi_1^0, \bvarphi_0^0)$ and $(\bvarphi_1^1,\bvarphi_0^1)$  are the solutions to \eqref{term-0} and  \eqref{term-1}, respectively.
\end{lem}

Recall that the domain $D$ is separated apart from  $ \Om$, then
\begin{align*}
\ds \sup_{x\in \Om, y\in \p D}\Big|\p^i \G_0(x-y)\Big|\leq C,\q i\in \NN^2,
\end{align*}
for some constant $C$ depending on $dist(D,\Om)$.  After the change of variables $\ye=\Phi_\ep(y)$,
we get from \eqref{sigexp}, \eqref{estimate-important}, and the Taylor expansion of $\G_0(x-\ye)$ in $y\in \p D$ for each fixed $x\in  \Om$ that
\begin{align*}
\ds \BScal_{0,D_\ep}[\widetilde{\bvarphi}_0](x)=&\int_{\p D_\ep}\G_0(x-\tilde y) \widetilde{\bvarphi}_0(\tilde y)d\sigma(\tilde y)\nonumber\\
\nm\ds =&\int_{\p D}\Big(\G_0(x-y)+\ep h(y)\nabla \G_0(x-y)\n(y)\Big)\Big( \bvarphi_0^0(y)+\ep \bvarphi^1_0(y)\Big)\nonumber\\
\nm\ds &\q\q\q\q\q\q\q\q\q\q\q\q\q \q\q \q\q \times \Big(1-\ep\kappa(y)h(y)\Big)d\sigma (y) +o(\ep)\nonumber\\
\nm\ds =& \BScal_{0,D}[\bvarphi_0^0](x)+\ep\Big(\BScal_{0,D}[\bvarphi_0^1](x)-\BScal_{0,D}[\kappa h\bvarphi_0^0](x)+\BDcal^{\sharp}_{0,D}[h\bvarphi_0^0](x)\Big)+o(\ep),
\end{align*}
Therefore, we obtain from  \eqref{representation-u-ep} that  for  $x\in \Om$,
 \begin{align}\label{main-equation-expansion}
\ds  \u_\ep(x) = \H(x)+ \BScal_{ 0, D}[\bvarphi^0_0](x) + \ep
 \Big(\BScal_{ 0, D}[\bvarphi^1_0](x) - \BScal_{ 0, D}[\kappa h \bvarphi^0_0](x)+ \BDcal^\sharp_{ 0, D} [h \bvarphi^0_0](x)\Big)+
 o(\ep).
 \end{align}
 According to \cite{book} (see also \cite{book2,AKNT}),   the solution $\u$ to \eqref{Main-Pb} is  represented as
\begin{equation}\label{equation-u-00}
\u(x)=\left\{
  \begin{array}{ll}
   \ds  \H(x)+\BScal_{0,D}[\bvarphi_0^0](x), &  x\in \RR^2\backslash \overline{D}, \\
    \nm \ds\BScal_{1,D}[\bvarphi_1^0](x),& x\in {D},
  \end{array}
\right.
\end{equation}
where $(\bvarphi_1^0,\bvarphi_0^0)$ is the unique solution of \eqref{term-0}.

The following theorem  follows immediately from \eqref{main-equation-expansion} and \eqref{equation-u-00}.
\begin{thm}\label{Thm-asymp-version0} For $\ep$ small enough. The following pointwise expansion holds  for $ x\in\Om$
\begin{align}\label{asymp-version0}
\ds \u_\ep(x)=\u(x)+\ep\Big(\BScal_{0,D}[\bvarphi^1_0](x)-\BScal_{0,D}[\kappa h\bvarphi_0^0](x)+\BDcal^{\sharp}_{0,D}[h\bvarphi_0^0](x)\Big)+o(\ep),
\end{align}
where $\bvarphi_0^0$ and $\bvarphi^1_0$ are defined by \eqref{term-0} and \eqref{term-1}, respectively.
The remainder $o(\ep)$ depends only  on $(\lambda_j, \mu_j)$ for j=0,1,2, the $\mathcal{C}^2$-norm of $X$,
the $\mathcal{C}^1$-norm of $h$, and $dist(\Om,D).$
\end{thm}

We now prove a representation theorem  for the solution of the transmission problem \eqref{equation-u-1} which will help us to derive the theorem \ref{Main-theorem}.
\begin{thm} \label{Representation-u1} The solution $\u_1$ of \eqref{equation-u-1} is represented by
\begin{equation}\label{u1-second}
\u_1(x)=
\left\{
  \begin{array}{lll}
   \ds  \BScal_{0,D}[\bvarphi^1_0](x)-\BScal_{0,D}[\kappa h\bvarphi_0^0 ](x)+\BDcal^\sharp_{0,D}[h\bvarphi_0^0](x),\q  x\in\RR^2 \backslash \overline{D}, \\
    \nm \ds \BScal_{1,D}[\bvarphi^1_1](x),\q  x\in D,  \\
  \end{array}
\right.
\end{equation}
where $\bvarphi^0_0$ and  $(\bvarphi^1_1,\bvarphi^1_0)$  are defined by  \eqref{term-0} and \eqref{term-1}, respectively.
\end{thm}
\proof  One can easily see that
\begin{align*}
\ds \Lcal_{\lambda_0,\mu_0}\u_1=0 \q \mbox{in }\RR^2\backslash \overline{D},\q\q\q \Lcal_{\lambda_1,\mu_1}\u_1=0\q \mbox{in }D.
\end{align*}
It follows from \eqref{term-1}, \eqref{equation-u-00}, and \eqref{identity-3} that
\begin{align*}
\ds \ds \u_1^i-\u_1^e=&\Big(\BScal_{1,D}[\bvarphi^1_1]-\BScal_{0,D}[\bvarphi^1_0]\Big)+\BScal_{0,D}[\kappa h\bvarphi^0_0 ]-\BDcal^\sharp_{0,D}[h\bvarphi_0^0]\big|_{+}\\
\nm\ds =&\Big[\mathcal{H}_1+\mathcal{Q}_1(\bvarphi_0^0)-\mathcal{Z}(\bvarphi_1^0)\Big]_{1}+\BScal_{0,D}[\kappa h\bvarphi^0_0 ]-\BDcal^\sharp_{0,D}[h\bvarphi_0^0]\big|_{+}\\
\nm\ds =& h\Big( \pd{\H}{\n}+\pd{\BScal_{0,D}[\bvarphi_{0}^0]}{\n}\Big|_{+}-\pd{\BScal_{1,D}[\bvarphi_{1}^0]}{\n}\Big|_{-}\Big)-h
\Big(\mathbb{K}_{2,1}\widehat{\nabla}\BScal_{0,D}[\bvarphi_{1}^0] \Big)\n\Big|_{-}\\
\nm\ds =& h\big(\nabla \u^e\n-\nabla\u^i \n\big)-h
\big(\mathbb{K}_{2,1}\widehat{\nabla}\u^i \big)\n\\
\nm\ds =&h \Big(\big(\mathbb{K}_{0,1}-\mathbb{K}_{2,1}\big)\widehat{\nabla} \u^i\Big)\n \q \mbox{on } \p D.
\end{align*}
Using \eqref{term-1} and \eqref{identity-1}, we get
\begin{align*}
\ds \ds \pd{\u_1}{{\nu_1}}\Big|_{-}-\pd{\u_1}{\nu_0}\Big|_{+}=&\Big(\pd{\BScal_{1,D}[\bvarphi^1_1]}{\nu_1}\Big|_{-}-
\pd{\BScal_{0,D}[\bvarphi^1_0]}{\nu_0}\Big|_{+}\Big)
+\pd{\BScal_{0,D}[\kappa h\bvarphi^0_0 ]}{\nu_0}\Big|_{+}-\pd{\BDcal^\sharp_{0,D}[h\bvarphi_0^0]}{\nu_0}\Big|_{+}\\
\nm\ds =&\Big[\mathcal{H}_1+\mathcal{Q}_1(\bvarphi_0^0)-\mathcal{Z}(\bvarphi_1^0)\Big]_{2}+\pd{\BScal_{0,D}[\kappa h\bvarphi^0_0 ]}{\nu_0}\Big|_{+}-\pd{\BDcal^\sharp_{0,D}[h\bvarphi_0^0]}{\nu_0}\Big|_{+}\\
\nm\ds =&\kappa h \pd{\H}{\nu_0}-\frac{\p}{\p \ta}\Big(h\big(\mathbb{C}_0\widehat{\nabla}\H)\big)\ta\Big)-\frac{\p}{\p \ta}\Big(h\big(\mathbb{C}_0\widehat{\nabla}\BScal_{0,D}[\bvarphi_0^0])\big)\ta\Big)\Big|_{+}\\
\nm\ds &+\kappa h \pd{\BScal_{0,D}[\bvarphi_0^0]}{\nu_0}\Big|_{+}-\kappa h \pd{\BScal_{1,D}[\bvarphi_1^0]}{\nu_1}\Big|_{-}+\frac{\p}{\p \ta}\Big(h\big(\mathbb{M}_{2,1}\widehat{\nabla}\BScal_{1,D}[\bvarphi_1^0])\big)\ta\Big)\Big|_{-}\\
\nm \ds =&\frac{\p}{\p \ta}\big(h (\mathbb{M}_{2,1}\widehat{\nabla} \u^i )\ta\big) -\frac{\p}{\p \ta}\big(h (\mathbb{C}_0\widehat{\nabla} \u^e )\ta\big)\\
\nm \ds=&\frac{\p}{\p \ta}\big(h ([\mathbb{M}_{2,1}-\mathbb{M}_{0,1}]\widehat{\nabla} \u^i )\ta\big).
\end{align*}
Now,  let us check the  condition
\begin{align}\label{phi1}
\ds \BScal_{0,D}[\bvarphi^1_0-\kappa h\bvarphi_0^0 ](x)\rightarrow 0\q \mbox{ as } |x|\rightarrow \infty.
\end{align}
To do this,  we rewrite the second component of the  equation  \eqref{term-1}
\begin{align}\label{phi-psi-ep-system-second}
 \ds \pd{\BScal_{1,D}[\bvarphi^1_1]}{ \nu_1}\Big |_{-}- \pd{{\BScal}_{0,D}[\bvarphi^1_0-\kappa h \bvarphi_0^0]}{ \nu_0}\Big|_{+}=&\frac{\p\BDcal_{0,D}^{\sharp}[h\bvarphi_0^0]}{\p  \nu_0}\Big|_{+}\nonumber \\
  \nm\ds &+\frac{\p}{\p\ta}\Big(h\big(\mathbb{M}_{2,1} \widehat{\nabla}\u^i\big)\ta\Big) -\frac{\p}{\p \ta}\Big(h\big(\mathbb{C}_0 \widehat{\nabla}\u^e\big)\ta\Big).
\end{align}
By \eqref{L-Psi},   $\p \BScal_{1,D}[\bvarphi^1_1]/\p{ \nu_1}\big |_{-} \in L^2_\Psi(\p D)$.  It is clear that
\begin{align*}
\ds \int_{\p D} \Big[\frac{\p}{\p\ta}\Big(h\big(\mathbb{M}_{2,1} \widehat{\nabla}\u^i\big)\ta\Big) -\frac{\p}{\p \ta}\Big(h\big(\mathbb{C}_0 \widehat{\nabla}\u^e\big)\ta\Big)\Big]\cdot \theta_m ~d\sigma =0\q \mbox{ for } m=1,2.
\end{align*}
Now,  we have
\begin{align}\label{007}
\ds \int_{\p D} \frac{d}{dt}\Big( h\big(\mathbb{C}_0 \widehat{\nabla}\u^e(x)\big)\ta(x)\Big) \cdot \theta_3(x) d\sigma& =-
\int_{\p D}h(x) \big(\mathbb{C}_0 \widehat{\nabla}\u^e(x)\big)\ta(x)\cdot \n(x) d\sigma\nonumber\\
\nm\ds &=-\mu_0\int_{\p D}h(x) \Big(\nabla \u^e(x)+(\nabla \u^e)^{T}(x) \Big)\ta(x)\cdot \n (x) d\sigma\nonumber\\
\nm\ds &=-\mu_0 \int_{\p D} h(x)\Big(\nabla \u^e(x)+(\nabla \u^e)^{T}(x) \Big)\n(x)\cdot \ta(x) d\sigma\nonumber\\
\nm\ds &=-\int_{\p D} h(x)\pd{\u^e}{\nu}(x)\cdot \ta (x)d\sigma.
\end{align}
We observe that
\begin{align}\label{observe}
\bvarphi_2=\bvarphi^0_2+O(\ep)\q \mbox{ and }\q \widetilde{\bpsi}_2\circ \Phi_\ep=\bpsi^0_2+O(\ep)\q \mbox{ on }\p D.
\end{align}
Substituting \eqref{observe} into the third line of \eqref{system}, we get
$$
 \Big(\mathbb{C}_2 \widehat{\nabla} \BScal_{2,D}[\bvarphi^0_2]\Big)\ta\Big |_{+}+\Big(\mathbb{C}_2 \widehat{\nabla} \BScal_{2,D}[{\bpsi}^0_2]\Big)\ta\Big |_{-}=(\mathbb{M}_{2,1} \widehat{\nabla}\u^i\Big)\Big |_{-}\ta.
$$
Similarly to \eqref{007}, one can easily see that
\begin{align*}
\ds \int_{\p D} \frac{\p}{\p \ta}\Big( h\big(\mathbb{M}_{2,1} \widehat{\nabla}\u^i\big)\ta \Big) \cdot \theta_3 d\sigma
=-\int_{\p D} h \Big( \pd{ \BScal_{2,D}[\bvarphi^0_2]}{\nu_2}\Big |_{+}+ \pd{\BScal_{2,D}[{\bpsi}^0_2]}{\nu_2}\Big |_{-}\Big)\cdot \ta d\sigma.
\end{align*}
By using ${\p\u_\ep}/{\p\nu_0}|_{+}={\p\u_\ep}/{\p\nu_2}|_{-}$ on $\p D_\ep$, we deduce from \eqref{asym-1}, \eqref{asym-5}, and \eqref{Expansion-conormal-H-xe}  that
$$
\pd{ \BScal_{2,D}[\bvarphi^0_2]}{\nu_2}\Big |_{+}+ \pd{\BScal_{2,D}[{\bpsi}^0_2]}{\nu_2}\Big |_{-}=\pd{\u^e}{\nu_0}\q \mbox{ on } \p D,
$$
which gives
\begin{align*}
\ds \int_{\p D} \Big[\frac{\p }{\p \ta }\Big(h\big(\mathbb{M}_{2,1} \widehat{\nabla}\u^i\big)\ta \Big) -\frac{\p}{\p \ta }\Big(h\big(\mathbb{C}_0 \widehat{\nabla}\u^e\big)\ta\Big)\Big]\cdot \theta_3d\sigma =0,
\end{align*}
and thus
\begin{align*}
\ds \frac{\p}{\p\ta}\Big(h\big(\mathbb{M}_{2,1} \widehat{\nabla}\u^i\big)\ta\Big) -\frac{\p}{\p \ta}\Big(h\big(\mathbb{C}_0 \widehat{\nabla}\u^e\big)\ta  \in L^2_\Psi(\p D).
\end{align*}
It then follows from \eqref{phi-psi-ep-system-second} that $\p \BScal_{0,D}[\bvarphi^1_0-\kappa h \bvarphi_0^0]/\p{ \nu_0}\big|_{+}+ \p\BDcal_{0,D}^{\sharp}[h\bvarphi_0^0]/\p \nu_0\big|_{+} \in L^2_\Psi(\p D)$. Since
\begin{align*}
\ds& \bvarphi^1_0-\kappa h \bvarphi_0^0+\frac{\p }{\p \ta}\Big(h ( \bvarphi_0^0\cdot \ta) \n +\frac{\lambda_0}{2\mu_0+\lambda_0}h (\bvarphi_0^0 \cdot \n ) \ta\Big)\\
\nm \ds &\q\q=
\pd{{\BScal}_D}{ \nu_0}\big[\bvarphi^1_0-\kappa h \bvarphi_0^0\big]\Big|_{+}+\frac{\p\BDcal_{0,D}^{\sharp}[h\bvarphi_0^0]}{\p \nu_0}\Big|_{+}-\pd{{\BScal}_D}{ \nu_0}\big[\bvarphi^1_0-\kappa h \bvarphi_0^0\big]\Big|_{-}-\frac{\p\BDcal_{0,D}^{\sharp}[h\bvarphi_0^0]}{\p \nu_0}\Big|_{-},
\end{align*}
with $\p \BScal_{0,D}[\bvarphi^1_0-\kappa h \bvarphi_0^0]/\p{ \nu_0}\big|_{-}, \p\BDcal_{0,D}^{\sharp}[h\bvarphi_0^0]/\p \nu_0\big|_{-} \in L^2_\Psi(\p D)$, see \eqref{L-Psi}. Then
\begin{align*}
\ds \bvarphi^1_0-\kappa h \bvarphi_0^0+\frac{\p }{\p \ta}\Big(h ( \bvarphi_0^0\cdot \ta) \n +\frac{\lambda_0}{2\mu_0+\lambda_0}h ( \bvarphi_0^0\cdot \n ) \ta\Big)\in L^2_\Psi(\p D).
\end{align*}
Therefore, we have
\begin{align*}
\ds \BScal_{0,D}[\bvarphi^1_0-\kappa h \bvarphi_0^0 ](x)&=\Gamma(x)\int_{\p D}(\bvarphi^1_0-\kappa h \bvarphi_0^0) d\sigma+O(|x|^{-1})\\
\nm \ds &=\Gamma(x)\int_{\p D}\bigg(\bvarphi^1_0-\kappa h \bvarphi_0^0+\frac{\p }{\p \ta}\Big(h ( \bvarphi_0^0 \cdot \ta) \n +\frac{\lambda_0}{2\mu_0+\lambda_0}h ( \bvarphi_0^0\cdot  \n ) \ta\Big)\bigg) d\sigma\\
\nm \ds&\q-\Gamma(x)\int_{\p D}\frac{\p }{\p \ta}\Big(h ( \bvarphi_0^0 \cdot \ta) \n +\frac{\lambda_0}{2\mu_0+\lambda_0}h ( \bvarphi_0^0\cdot  \n ) \ta\Big) d\sigma+O(|x|^{-1})\\
\nm\ds &=O(|x|^{-1})\q \mbox{ as }|x|\rightarrow \infty.
\end{align*}
Thus   $\u_1$ defined by   \eqref{u1-second} satisfies $\u_1(x)=O(|x|^{-1})$   as $|x|\rightarrow \infty.$ This completes the proof of the theorem \ref{Representation-u1}.

The main theorem \ref{Main-theorem}  immediately follows from  the integral representation of $\u_1$   in  \eqref{u1-second} and the theorem \ref{Thm-asymp-version0}.
\subsection{Proof of the theorem \ref{second-theorem}}
The following corollary can be proved as in exactly the same manner as Theorem \ref{Main-theorem}.
\begin{cor} Let $\u_\ep$ and $\u$ be the solutions to
\eqref{equation-u} and \eqref{Main-Pb}, respectively. Let $\Om$ be a bounded
region away from $\p D$. For $x \in \Om$, the following  pointwise asymptotic expansion holds:
\begin{equation}\label{Main-Asymptotic-normal}
\pd{\u_\ep}{\nu_0}(x)=\pd{\u}{\nu_0}(x)+\ep \pd{\u_1}{\nu_0}(x)+o(\ep),
\end{equation}
where the remainder $o(\ep)$ depends only  on $(\lambda_j, \mu_j)$ for j=0,1,2,
the $\mathcal{C}^2$-norm of $X$, the $\mathcal{C}^1$-norm of $h$, and $dist (\Om, \p D)$,
and  $\u_1$ is the unique solution of \eqref{equation-u-1}.
\end{cor}
 Let us note  simple, but important relations.
\begin{lem}\label{relations}
\begin{enumerate}
 \item  If $\f\in W^{1,2}(D)$ and $\Lcal_{\lambda_0,\mu_0}\f=0$ in $D$, then for all $\g \in W^{1,2}(D)$,
\begin{align}\label{Important-relation}
 \ds \displaystyle  \int_{\p D} \g\cdot \pd{\f}{{\nu_0}} ~d\sigma=\int_{D} \lambda_0 (\nabla \cdot \f)(\nabla \cdot \g)+\frac{\mu_0}{2}(\nabla \f+\nabla \f ^{T}): (\nabla \g+\nabla \g ^{T}) d\sigma.
\end{align}
 \item If $\f \in W^{1,2}(\RR^2\backslash \overline{D})$ and $\Lcal_{\lambda_0,\mu_0}\f=0$ in $\RR^2\backslash \overline{D}$, $\f(x)=O(|x|^{-1})$  as $|x|\rightarrow \infty$. Then for all $\g \in W^{1,2}(\RR^2\backslash \overline{D})$,  $\g(x)=O(|x|^{-1})$  as $|x|\rightarrow \infty$, we have
\begin{align}\label{Important-relation-2}
\ds  - \displaystyle  \int_{\p D} \g\cdot \pd{\f}{{\nu_0}} ~d\sigma=\int_{\RR^2\backslash \overline{D}} \lambda_0 (\nabla \cdot \f)(\nabla \cdot \g)+\frac{\mu_0}{2}(\nabla \f+\nabla \f ^{T}): (\nabla \g+\nabla \g ^{T}) d\sigma.
\end{align}
Here,  for $2\times 2$ matrices $\M$ and $\N$,  $\M: \N= \displaystyle \sum_{ij}\M_{ij}\N_{ij}$.
\end{enumerate}
\end{lem}

 Let $S$ be a Lipschitz closed curve enclosing $D$  away from $\p D$. Let $\v$ be the solution to \eqref{v}. It follows from \eqref{Main-asymptotic}, \eqref{Main-Asymptotic-normal}, and \eqref{Important-relation-2} that
\begin{align*}
\ds\int_{S}\big(\u_\ep-\u\big)\cdot\pd{ \F}{\nu_0}d\sigma-\int_{S}\big(\pd{\u_\ep}{\nu_0}-
\pd{ \u}{\nu_0}\big) \cdot\F d\sigma& =\ep \int_{S}\Big(\u_1\cdot \pd{ \v}{\nu_0}-\pd{\u_1}{\nu_0}
\cdot \v\Big) d\sigma+o(\ep).
\end{align*}
By using Lemma \ref{relations} to  the integral on the right-hand side, we get
$$
\int_{S}\Big(\u_1\cdot \pd{ \v}{\nu_0}-\pd{\u_1}{\nu_0}
\cdot \v\Big) d\sigma=
 \int_{\p D} \Big(\pd{\v^{e}}{\nu_0}\cdot \u^{e}_1-\v^e \cdot\pd{\u^{e}_1}{\nu_0}\Big)d\sigma.
$$
According to the jump conditions for $\u_1$ in \eqref{equation-u-1}, we deduce that
\begin{align}\label{eq001}
\ds\int_{S}\Big(\u_1\cdot \pd{ \v}{\nu_0}-\pd{\u_1}{\nu_0}
\cdot \v\Big) d\sigma=& \int_{\p D} \Big(\pd{\v^{i}}{\nu_1}\cdot \u^{i}_1-\v^i \cdot \pd{\u^{i}_1}{\nu_1}\Big)d\sigma\nonumber\\
\nm\ds &-\int_{\p D} h\Big(\big(\mathbb{K}_{0,1}-\mathbb{K}_{2,1}\big)\widehat{\nabla}\u^{i}\Big)\n \cdot \big(\mathbb{C}_1\widehat{\nabla}\v^{i}\big)\n d\sigma\nonumber\\
\nm \ds &+ \int_{\p D} \frac{\p}{\p \ta }\Big(h\big([\mathbb{M}_{2,1}-\mathbb{M}_{0,1}]\widehat{\nabla}\u^{i}\big)\ta \Big) \cdot \v^{i} d\sigma .
\end{align}
It follows from \eqref{Important-relation} that
\begin{align}\label{eq002}
\ds\int_{\p D} \Big(\pd{\v^{i}}{\nu_1}\cdot \u^{i}_1-\v^i \cdot \pd{\u^{i}_1}{\nu_1}\Big)d\sigma=0.
\end{align}
We have
\begin{align}\label{eq003}
 \ds \int_{\p D} \frac{\p}{\p \ta } \Big(h\big([\mathbb{M}_{2,1}-\mathbb{M}_{0,1}]\widehat{\nabla}\u^{i}\big)\ta \Big) \cdot \v^{i} d\sigma &=
-\int_{\p D}  h\big([\mathbb{M}_{2,1}-\mathbb{M}_{0,1}]\widehat{\nabla}\u^{i}\big)\ta  \cdot \nabla \v^{i} \ta d\sigma.
\end{align}
One can easily check that
\begin{align}\label{eq004}
 \ds \big([\mathbb{M}_{2,1}-\mathbb{M}_{0,1}]\widehat{\nabla}\u^{i}\big)\ta  \cdot \nabla \v^{i} \ta
=\big([\mathbb{M}_{2,1}-\mathbb{M}_{0,1}]\widehat{\nabla}\u^{i}\big)\ta  \cdot \widehat{\nabla} \v^{i} \ta.
\end{align}
We finally obtain from \eqref{eq001}-\eqref{eq004}
the relationship between traction-displacement measurements and the shape deformation $h$ \eqref{asymptotic-traction-displacement}, as desired.

\newpage
\section*{Appendix }
 Let  $\w$ be   the solution of  $ \nabla \cdot \big(\mathbb{C} \widehat{\nabla}\w \big)=0$  in  $\RR^2$, where $\mathbb{C}:=\mathbb{C}_l\chi_{\RR^2\backslash \overline{D}}+\mathbb{C}_k\chi_{D}.$  Then $\w$ satisfies the transmission
conditions along the interface $\p D$:
\begin{align}
\ds \ds  \w^i&=\w^e,\\
\nm  \ds \nabla \w^i \ta &= \nabla \w^e \ta, \label{eq01} \\
  \nm \ds \la \widehat{\nabla} \w^i \ta, \ta\ra  &=   \la \widehat{\nabla} \w^e \ta, \ta\ra,\label{eq02} \\
\nm \ds \lambda_k \nabla \cdot \w^i+2 \mu_k  \la \widehat{\nabla} \w^i \n, \n\ra&=
\lambda_l \nabla \cdot \w^e+2 \mu_l  \la \widehat{\nabla} \w^e \n, \n\ra,\label{eq03}\\
\nm \ds \mu_k  \la \widehat{\nabla} \w^i \n, \ta\ra&=\mu_l  \la \widehat{\nabla} \w^e \n, \ta\ra.\label{eq04}
\end{align}
We have from \cite{JFH} the following lemma.
\begin{lem} \label{Important-prop} We have the following identities  along the interface $\p D$
\begin{align}
 \ds \big(\mathbb{C}_l \widehat{\nabla} \w^e\big)\ta&=\big(\mathbb{M}_{l,k} \widehat{\nabla} \w^i\big)\ta,\label{identity-1}\\
 \nm\ds \big(\mathbb{C}_k \widehat{\nabla} \w^i\big)\ta&=\big(\mathbb{M}_{k,l} \widehat{\nabla} \w^e\big)\ta,\label{identity-2}\\
\nm\ds \nabla \w^e \n -\nabla \w^i \n&=\big(\mathbb{K}_{l,k} \widehat{\nabla} \w^i\big)\n=-\big(\mathbb{K}_{k,l} \widehat{\nabla} \w^e\big)\n,\label{identity-3}
\end{align}
where the  $4$-tensors $\mathbb{M}_{l,k}$ and $\mathbb{K}_{l,k}$ are defined by:
\begin{align*}
\ds \mathbb{M}_{l,k}&:=\frac{\lambda_l (\lambda_k+2\mu_k)}{\lambda_l+2\mu_l} \I \otimes \I+2 \mu_k \mathbb{I}+ \frac{4(\mu_l-\mu_k)(\lambda_l+\mu_l)}{\lambda_l+2\mu_l} \I\otimes (\ta\otimes \ta),\\
\nm \ds \mathbb{K}_{l,k}&:=\frac{\mu_l(\lambda_k-\lambda_l)+2(\mu_l-\mu_k) (\lambda_l+\mu_l)}{\mu_l(\lambda_l+2\mu_l)} \I \otimes \I+2\big(\frac{\mu_k}{\mu_l}-1\big) \mathbb{I}\\
\nm\ds &\q+ \frac{2(\mu_k-\mu_l)(\lambda_l+\mu_l)}{\mu_l(\lambda_l+2\mu_l)} \I\otimes (\ta\otimes \ta).
\end{align*}
\end{lem}
\proof  Recalling that
\begin{align}\label{Tr-w}
\ds \nabla \cdot \w^e=\widehat{\nabla} \w^{e}: \I= tr(\widehat{\nabla} \w^{e})=\la\widehat{\nabla} \w^{e}\n,\n \ra+ \la\widehat{\nabla} \w^{e}\ta,\ta \ra.
\end{align}
Here $\la, \ra$ denotes the scalar product in $\RR^2$.

One can easily get from \eqref{eq02}, \eqref{eq03}, and \eqref{Tr-w} that
\begin{align}\label{nabla-cdot-w}
\ds \nabla \cdot \w^e= \frac{\lambda_k+2\mu_k}{\lambda_l+2\mu_l}\nabla \cdot \w^i+\frac{2(\mu_l-\mu_k)}{\lambda_l+2\mu_l}\la \widehat{\nabla} \w^i \ta ,\ta\ra.
\end{align}
We have
\begin{align*}
\ds \nabla \w^{e} \n&=\la\nabla \w^{e} \n, \n\ra \n +\la\nabla \w^{e} \n, \ta\ra \ta\\
\nm\ds &=\la\widehat{\nabla} \w^{e} \n, \n\ra \n +2\la\widehat{\nabla} \w^{e} \n, \ta\ra \ta-\la(\nabla \w^{e})^{T} \n, \ta\ra \ta\\
\nm\ds &=\la\widehat{\nabla} \w^{e} \n, \n\ra \n +2\la\widehat{\nabla} \w^{e} \n, \ta\ra \ta-\la \nabla \w^{e} \ta, \n\ra \ta.
\end{align*}
According to \eqref{Tr-w}, we get
\begin{align*}
\ds \nabla \w^{e} \n=(\nabla \cdot \w^{e}) \n -\la\widehat{\nabla} \w^{e} \ta, \ta\ra \n +2\la\widehat{\nabla} \w^{e} \n, \ta\ra \ta-\la \nabla \w^{e} \ta, \n\ra \ta.
\end{align*}
In a similar way as before, we write
\begin{align*}
\ds \nabla \w^{i} \n&=(\nabla \cdot \w^{i})\n -\la\widehat{\nabla} \w^{i} \ta, \ta\ra \n +2\la\widehat{\nabla} \w^{i} \n, \ta\ra \ta-\la \nabla \w^{i} \ta, \n\ra \ta.
\end{align*}
It then follows from \eqref{eq04}, \eqref{Tr-w}, and \eqref{nabla-cdot-w} that
\begin{align*}
\ds \nabla \w^{ e} \n-\nabla \w^{i}\n&=(\nabla \cdot \w^{e}-\nabla \cdot \w^{i})\n+2\la\widehat{\nabla} \w^{e} \n, \ta\ra \ta-2\la\widehat{\nabla} \w^{i} \n, \ta\ra \ta\\
\nm \ds &=\Big(\frac{\lambda_k+2\mu_k}{\lambda_l+2\mu_l}-1\Big)(\nabla \cdot \w^{i})\n+2\big(\frac{\mu_k}{\mu_l}-1\big)\la\widehat{\nabla} \w^{i} \n, \ta\ra \ta \\
\nm\ds &\q+\frac{2(\mu_l-\mu_k)}{\lambda_l+2\mu_l} \la \widehat{\nabla} \w^i \ta, \ta\ra \n\\
\nm \ds &=\Big(\frac{\lambda_k+2\mu_k}{\lambda_l+2\mu_l}-1-2\big(\frac{\mu_k}{\mu_l}-1\big)\Big)(\nabla \cdot \w^{i})\n+2\big(\frac{\mu_k}{\mu_l}-1\big)\widehat{\nabla} \w^{i} \n \\
\nm\ds &\q+\Big(\frac{2(\mu_l-\mu_k)}{\lambda_l+2\mu_l}+2\big(\frac{\mu_k}{\mu_l}-1\big)\Big) \la \widehat{\nabla} \w^i \ta, \ta\ra \n\\
\nm \ds &=\frac{\mu_l(\lambda_k-\lambda_l)+2(\mu_l-\mu_k) (\lambda_l+\mu_l)}{\mu_l(\lambda_l+2\mu_l)}(\nabla \cdot \w^{i})\n+2\big(\frac{\mu_k}{\mu_l}-1\big)\widehat{\nabla} \w^{i} \n \\
\nm\ds &\q+ \frac{2(\mu_k-\mu_l)(\lambda_l+\mu_l)}{\mu_l(\lambda_l+2\mu_l)} \la \widehat{\nabla} \w^i \ta, \ta\ra \n\\
\nm \ds &=\big(\mathbb{K}_{l,k} \widehat{\nabla} \w^i\big)\n\q\mbox{on }\p D.
\end{align*}
We obtain from  \eqref{eq02}-\eqref{eq04}, and \eqref{nabla-cdot-w} that
\begin{align*}
\ds   \big(\mathbb{C}_l \widehat{\nabla} \w^e \big)\ta&= \lambda_l (\nabla \cdot \w^e )\ta+2\mu_l (\widehat{\nabla} \w^e) \ta\\
\nm \ds &=\frac{\lambda_l (\lambda_k+2\mu_k)}{\lambda_l+2\mu_l}(\nabla \cdot \w^i)\ta+\frac{2\lambda_l(\mu_l-\mu_k)}{\lambda_l+2\mu_l}\la \widehat{\nabla} \w^i \ta, \ta\ra \ta\\
\nm \ds &\q + 2\mu_l\la \widehat{\nabla} \w^i \ta, \ta\ra \ta+2\mu_k\la \widehat{\nabla} \w^i \ta, \n\ra \n\\
\nm \ds &=\frac{\lambda_l (\lambda_k+2\mu_k)}{\lambda_l+2\mu_l}(\nabla \cdot \w^i)\ta+\frac{2\lambda_l(\mu_l-\mu_k)}{\lambda_l+2\mu_l}\la \widehat{\nabla} \w^i \ta, \ta\ra \ta\\
\nm \ds &\q +2\mu_k \widehat{\nabla }\w^i \ta + 2\mu_l\la \widehat{\nabla} \w^i \ta, \ta\ra \ta- 2\mu_k\la \widehat{\nabla} \w^i \ta, \ta\ra \ta\\
\nm \ds &= \frac{\lambda_l (\lambda_k+2\mu_k)}{\lambda_l+2\mu_l} (\nabla \cdot \w^i)\ta+2\mu_k \widehat{\nabla }\w^i \ta +\frac{4(\mu_l-\mu_k)(\lambda_l+\mu_l)}{\lambda_l+2\mu_l} \la \widehat{\nabla} \w^i \ta, \ta\ra \ta\\
\nm\ds &=\big(\mathbb{M}_{l,k} \widehat{\nabla} \w^i \big)\ta\q\mbox{on }\p D.
\end{align*}
The identities $(\mathbb{C}_k \widehat{\nabla} \w^i\big)\ta=\big(\mathbb{M}_{k,l} \widehat{\nabla} \w^e\big)\ta$ and $\nabla \w^e \n -\nabla \w^i \n=-\big(\mathbb{K}_{k,l} \widehat{\nabla} \w^e\big)\n$ can be done in exactly the same manner as above.

\end{document}